\documentclass[11pt]{article}
\usepackage{amssymb}
\usepackage{latexsym}

\newtheorem{thm}{Theorem}[section]
\newtheorem{lem}[thm]{Lemma}
\newtheorem{cor}[thm]{Corollary}
\newtheorem{pro}[thm]{Proposition}
\newtheorem{ex}[thm]{Example}

\newtheorem{defi}[thm]{Definition}

\newcommand{\gm }{\Gamma }
\newcommand{\lon }{\longrightarrow }
\newcommand{\be }{\begin{eqnarray*}}
\newcommand{\ee }{\end{eqnarray*}}

\newcommand{\per }{\backl }

\input amssym.def
\input amssym

\newcommand{\pf}{\noindent{\bf Proof.}\ }
\newcommand{\qed}{\begin{flushright} $\Box$\ \ \ \ \ \
                  \end{flushright}}
\newcommand{\complex}{{\Bbb C}}
\newcommand{\reals}{{\Bbb R}}

\newcommand{\frakg}{{\frak g}}

\newcommand{\half}{\frac{1}{2}}

\newcommand{\backl}{\mathbin{\vrule width1.5ex height.4pt\vrule height1.5ex}}
\newcommand{\cala}{{\cal A}}

\newcommand{\cald}{{\cal D}}

\newcommand{\calf}{{\cal F}}

\newcommand{\calp}{{\cal P}}

\newcommand{\calx}{{\cal X}}

\newcommand{\smalcirc}{\mbox{\tiny{$\circ $}}}

\def\description label#1{\hfil\bf[#1]\hfil}
\newcommand{\parr}[1]{\frac{\partial}{\partial x^{#1}}}
\newcommand{\parrr}[2]{\frac{\partial #1}{\partial x^{#2}}}

\newcommand{\ih}{\frac{i}{\hbar}}
\newcommand{\poll}{pol (\frakg^{*})}
\newcommand{\pol}{U\frakg }
\newcommand{\ad}{ad}
\newcommand{\muh}{\mu_{\hbar}}
\newcommand{\commutant}{commutant }
\newcommand{\equivariant}{equivariant }
\newcommand{\preferred}{regular }
\newcommand{\nat}{natural }

\def\sdp{\mathbin{\hbox{$\mapstochar\kern-.3333em\times$}}}
\def\pds{\mathbin{\hbox{$\times\kern-.55em\mapstochar\,$}}}

\newcommand{\wed}{\mathbin{\lower1.5pt\hbox{$\scriptstyle{\wedge}$}}}

\let\Tilde=\widetilde

\def\chigh{{\raise1.5pt\hbox{$\chi$}}}
\let\phi=\varphi
\def\til0{\Tilde{0}}

\def\dminus{\raise2pt\hbox{\vrule height1pt width 2ex}\hskip3pt}

\def\pback#1{\mathbin{{{\lower1.2ex\hbox{$\times$}}\atop #1}}}

\def\vlra{\hbox{$\,-\!\!\!-\!\!\!-\!\!\!-\!\!\!-\!\!\!
-\!\!\!-\!\!\!-\!\!\!-\!\!\!-\!\!\!\longrightarrow\,$}}

\def\gpd{\,\lower1pt\hbox{$\longrightarrow$}\hskip-.24in\raise2pt
             \hbox{$\longrightarrow$}\,}

\def\lgpd{\,\lower1pt\hbox{$\vlra$}\hskip-1.02in\raise2pt\hbox{$\vlra$}\,}

\def\llgpd{\,\lower1pt\hbox{$\vvlra$}\hskip-1.3in\raise2pt\hbox{$\vvlra$}\,}

\begin{document}

\title{{\bf  Fedosov $*$-products and quantum momentum maps}}

\author{ PING XU \thanks{ Research partially supported by NSF
        grant DMS95-04913}\\
 Department of Mathematics\\The  Pennsylvania State University \\
University Park, PA 16802, USA\\
        {\sf email: ping@math.psu.edu }}

\date{August, 1996}

\maketitle

\begin{abstract}
The purpose of the  paper  is to study various aspects of  star products on
a symplectic manifold related to the  Fedosov method.

 By introducing the notion of     ``quantum exponential maps",
 we give  a criterion    characterizing Fedosov connections. 
As a consequence, a geometric realization is
obtained  for the equivalence between an arbitrary $*$-product and
a Fedosov one.

 Every Fedosov $*$-product is  shown to be  a Vey $*$-product.
Consequently, one obtains that
 every $*$-product is equivalent to a Vey $*$-product,  a classical result
of  Lichnerowicz.

Quantization of a hamiltonian $G$-space, and in particular, quantum
momentum maps are studied. Lagrangian submanifolds are  also
studied under a deformation quantization.
\end{abstract}

\section{Introduction}

In classical mechanics, observables are  smooth functions on  a phase
space,  which consist of a Poisson algebra, while in quantum mechanics
observables become a noncommutative associative algebra.
Deformation quantization, as laid out by
Bayen, Flato, Fronsdal, Lichnerowicz and Sternheimer in 1970's \cite{BFFLS}, is
one of the important  attempts
 aiming to establish a correspondence  principle between these two mechanics.
A classical phase space is usually  a symplectic  manifold $M$.
A deformation quantization, or more precisely a star-product is a
family of associative  multiplication $*_{\hbar}$ (depending on the Planck
constant $\hbar$) on $C^{\infty}(M)[[\hbar]]$, the
space of formal power series with coefficients in $C^{\infty}(M)$:
\begin{equation}
f*_{\hbar}g=fg  -\frac{i\hbar }{2}\{f , g\} + \cdots +\hbar^{k}C_{k}(f, g)+\cdots  , \ \ \forall
f, g\in C^{\infty}(M)[[\hbar]]
\end{equation}
such that
\begin{enumerate}
\item $C_{k}(f, g)=(-1)^{k}C_{k}(g, f)$;
\item  $C_{k}(1, f)=C_{k}(f, 1)=0$, for $k\geq 1$;
\item  each $C_{k}(\cdot , \cdot )$ is a  bidifferential operator.
\end{enumerate}
Here  $\{f , g\}$ is the Poisson bracket on the symplectic manifold $M$.

Such a $*$-product always exists on   a symplectic vector space $V$, which is
known as  Moyal-Weyl formula:
\begin{equation}
\label{eq:mw}
f  *_{\hbar} g=\sum_{k=0}^{\infty} (-\frac{i\hbar}{2})^{k}\frac{1}{k!}\pi^{i_{1}j_{1}}
\cdots \pi^{i_{k}j_{k}}
\frac{\partial^{k} f}{\partial y^{i_{1}}\cdots \partial y^{i_{k}}}
\frac{\partial^{k} g}{\partial y^{j_{1}}\cdots \partial y^{j_{k}}}, \ \ \forall f , g\in C^{\infty}(V)[[\hbar ]],
\end{equation}
where $y^{1}, \cdots , y^{2n}$ are linear coordinates on $V$, and $\pi^{ij}=\{y^{i} , y^{j}\}$.
It is simple to see that this definition is independent of
the choice  of linear coordinates.

Moyal-Weyl formula has  a straightforward generalization
to the case of a symplectic manifold admitting  a flat torsion free
symplectic connection $\nabla$, as shown in \cite{BFFLS}.
In this case, for any $f, \ g\in C^{\infty}(M)[[\hbar ]]$, set
\begin{equation}
\label{eq:generalized-mw}
(f*_{\hbar }g)(x)=[(exp^{*}_{x}f)(y)*_{\hbar }(exp^{*}_{x}f)(y)]|_{y=0},
\end{equation}
where $exp_{x}: T_{x}M\lon M$ is the exponential map, defined
in a neighborhood of the origin, corresponding to the connection
$\nabla$, and the $*_{\hbar }$-product on the RHS refers to the standard
Moyal-Weyl $*$-product on the symplectic vector space $T_{x}M$.
More explicitly, one has
\begin{equation}
(f*_{\hbar }g)(x) =
\sum_{k=0}^{\infty} (-\frac{i\hbar}{2})^{k}\frac{1}{k!}
\pi^{i_{1} j_{1}}\cdots \pi^{ i_{k} j_{k}}(\partial_{i_{1}}\cdots \partial_{i_{k}}f)
(\partial_{j_{1}}\cdots \partial_{j_{k}}g).
\end{equation}

However, the multiplication defined by Equation (\ref{eq:generalized-mw})
fails to be associative when $\nabla$ has  a curvature.
The  existence proof  of $*$-products on a general
symplectic manifold was first obtained by de Wilde and Lecomte \cite{DeL}
using  a   homological argument. Later, an alternate proof using Weyl
manifolds was found by Omori  et al.\cite{OMY}.
Guillemin showed that their results in \cite{GV} also
implies the existence of $*$-products.

Recently, Fedosov has given a nice geometrical   existence proof \cite{F1}
\cite{F2}, which
provides a useful tool  for understanding $*$-products geometrically.
This paper grew out from an attempt to understand the Fedosov method, as well as
an attempt to solve   some other  related  problems
  using such a  method.

Roughly speaking, Fedosov's method is to make a ``quantum correction"
to Equation  (\ref{eq:generalized-mw}) when the   curvature exists.
As discussed early in this  introduction, each tangent space $T_{x}M$ is
a symplectic vector space, hence can be quantized  using the standard
Moyal-Weyl product. This enables us to obtain a bundle of algebras $W\lon M$,
called Weyl bundle by Fedosov, which can be thought of  as
a kind of ``quantum tangent bundle". Fedosov found a nice iteration method
 of constructing a flat connection on the Weyl bundle, whose
parallel sections can be naturally identified with
 $C^{\infty}(M)[[\hbar ]]$. Thus the product on the bundle induces
a $*$-product on $C^{\infty}(M)[[\hbar ]]$. Intuitively, Fedosov
connection  can be thought of as a ``quantum connection" on  the ``quantum
tangent bundle" $W$,
which is indeed obtained by adding some ``quantum correction" to
the usual affine connection on the tangent bundle (see Section 2).
Then,  the correspondence between $C^{\infty}(M)[[\hbar ]]$ and
its parallel sections can be considered as taking  the exponential map for such
a quantum connection. This viewpoint  has already been  pointed out by Weinstein
(see \cite{Weinstein:1995} \cite{EW}), and is essential   to our introduction
of the  notion of  ``quantum exponential maps".

By a ``quantum exponential map", we mean a map from $C^{\infty}(M)[[\hbar ]]$
to the space of sections $\gm W$ of the Weyl bundle, which
satisfies certain axioms (see Section 3 for details). In fact, we
will show that any quantum exponential map is  equivalent
to a Fedosov connection, which is what is expected: connections
and exponential maps are equivalent.  This result will  enable us to
 characterize  those    subalgebras  of $\gm W$ arising from
 Fedosov connections. As an application, we will 
 realize    geometrically
 the equivalence between any  given $*$-product and
a Fedosov  one, a result  proved  by several authors including
 Nest-Tsygan \cite{NT1},  Deligne \cite{Deligne} and recently Bertelson-Cahen-Gutt \cite{BCG}.   

Since  a ``quantum exponential map" is a  ``quantum correction" to the usual
exponential map, the $*$-products obtained by Fedosov method, called
Fedosov $*$-products in the paper,
 should be closely related to
 Equation  (\ref{eq:generalized-mw}). A Vey $*$-product,
by definition,  is a $*$-product
where the principal terms coincide with those as in Equation (\ref{eq:generalized-mw})
(see Section 4 for the precise definition).
Vey $*$-products have played an important role since the beginning of 
deformation quantization  theory (see \cite{BFFLS} \cite{L} \cite{V}).
In this paper, we will show that every Fedosov $*$-product
is a Vey $*$-product. As a consequence, we recover
the following well known result of  Lichnerowicz \cite{L}:
 any $*$-product is equivalent to a Vey $*$-product.

Because of the  basic simplicity of its construction,
    Fedosov method provides us a useful tool for  studying some other problems
in deformation quantization  theory.  In this paper,
in particular, we will study the following question:

\begin{quote}
What do lagrangian submanifolds correspond under a deformation
quantization?
\end{quote}

Given a lagrangian submanifold $L$, the space  $C^{\infty}_{L}(M)$ of  smooth functions
vanishing on $L$ forms a Poisson subalgebra. For any
 $*$-product,  we will show that under some ``quantum
correction" $C^{\infty}_{L}(M)[[\hbar ]]$ becomes a subalgebra.

In symplectic geometry, hamiltonian $G$-spaces, and in particular, momentum
maps play a very special role. It is natural to ask: what
is the quantum analogue of momentum maps? The second part of the paper,
as another application of Fedosov method,
is devoted to the study of  quantum momentum maps.
In particular, we will derive a sufficient condition
for their  existence.
When a quantum momentum map exists, we obtain a  pair of
mutual commutants  which can be considered as a 
quantum analogue of the well known Poisson  dual pair $\frakg^{*}\stackrel{J}{\longleftarrow}M
\stackrel{pr}{\lon}M/G$ of Weisntein \cite{Weinstein:1983}.

Derivations are always important for  associative algebras. In Appendix A,
we will collect some basic  results regarding
 derivations on a $*$-algebra $C^{\infty}(M)[[\hbar ]]$,
which are needed for the study of quantum momentum maps.
Most of them  are proved  using Fedosov method. Even for
some well known results, we will see that Fedosov method provides 
a simpler way in understanding them.

While the paper is in writing,
it came to the author's attention that quantum momentum  maps are also being
  under study by some other authors including Astashkevich \cite{Ast}, Kostant
and Tsygan. In particular,  Theorem \ref{eq:mom-exist} has  also been  proved by Kostant. 

{\bf Acknowledgements.} The author  would like to thank Pierre Bieliavsky, Jean-Luc Brylinski,
Ranee Brylinski,  Moshe Flato, 
Victor Guillemin, Yuri Manin, Marc Rieffel, Daniel Sternheimer,
Boris Tsygan and Alan Weinstein for useful discussions.  In
addition to the funding sources mentioned
in the first footnote, he would also  like to thank IHES and Max-Planck-Institut
for their hospitality and financial support while part of this project was being done.

\section{Fedosov quantization}
In this section, we will recall some basic
ingredients  of
Fedosov construction of
$*$-products on a symplectic manifold, as well as some
useful notations,  which will
be needed in the future. For details, readers should
consult \cite{F1} \cite{F2}.

Let $(M, \omega )$ be a symplectic manifold of dimension $2n$.
Then, each tangent space $T_{x}M$ is equipped with  a linear symplectic
structure, which can be quantized by the standard Moyal-Weyl product.
The resulting space is denoted by $W_{x}$.
More precisely,

\begin{defi}
A formal Weyl  algebra $W_{x}$  associated to $T_{x}M$
is an associative algebra  with a unit over $\complex $,
whose elements consist of formal power series in the formal parameter
$\hbar $ with coefficients being formal  polynomials
in $T_{x}M$. In other words,
each element has the form
\begin{equation}
\label{eq:general}
a(y, \hbar )=\sum \hbar^{k}a_{k, \alpha }y^{\alpha }
\end{equation}
where   $y=(y^{1}, \cdots , y^{2n})$ is a linear coordinate
on $T_{x}M$, $\alpha =(\alpha_{1}, \cdots , \alpha_{2n})$
is a multi-index and         $y^{\alpha}=(y^{1})^{\alpha_{1}}\cdots
(y^{2n})^{\alpha_{2n}}$.
The product is defined by the Moyal-Weyl rule:

\begin{equation}
\label{eq:moyal}
a  *b=\sum_{k=0}^{\infty} (-\frac{i\hbar}{2})^{k}\frac{1}{k!}\pi^{i_{1}j_{1}}
\cdots \pi^{i_{k}j_{k}}
\frac{\partial^{k} a}{\partial y^{i_{1}}\cdots \partial y^{i_{k}}}
\frac{\partial^{k} b}{\partial y^{j_{1}}\cdots \partial y^{j_{k}}} .
\end{equation}
\end{defi}


Let $W=\cup_{x\in M}W_{x}$. Then $W$ is a bundle of algebras over $M$, called
the Weyl bundle over $M$.
 Its  space of   sections $\gm W$  forms  an associative algebra
with unit under
the fiberwise multiplication.
One may think of $W$ as a ``quantum tangent bundle"
of  $M$, whose space of sections $\gm W$ gives
rise to a deformation quantization for the
tangent bundle $TM$ considered as a Poisson manifold,
with fiberwise linear symplectic structure on  $T_{x}M, \ \ x\in M$.

 The center $Z (W)$  of $\gm W$ consists of sections
 not containing $y's$, thus can be naturally
identified with $C^{\infty}(M)[[\hbar ]]$.

By assigning degrees to $y's$ and $\hbar$
with $deg y^{i}=1$ and $deg \hbar =2$, there
is a natural filtration
$$\gm (W)\subset \gm (W_{1})\subset \cdots \gm (W_{i})\subset \gm (W_{i+1})\cdots   $$ with respect
to the total degree $2k+l$ of the series terms  in Equation (\ref{eq:general}).

A differential form with values in $W$ is a section of the
bundle $W\otimes \wedge^{q}T^{*}M$, which  
can be expressed locally as
\begin{equation}
a(x, y, \hbar , dx)=\sum \hbar^{k}   a_{k, i_{1}\cdots i_{p}, j_{1}\cdots j_{q}}
y^{i_{1}}\cdots y^{i_{p}}dx^{j_{1}}\wedge \cdots \wedge dx^{j_{q}}.
\end{equation}
Here the coefficient $a_{k, i_{1}\cdots i_{p}, j_{1}\cdots j_{q}} $
is a covariant tensor symmetric with respect to
$i_{1} \cdots i_{p}$ and antisymmetric in $j_{1}\cdots j_{q}$.
For short,  we denote the space of
sections of the bundle $W\otimes \wedge^{q}T^{*}M$ by
$\gm W\otimes \Lambda^{q}$.

The usual   exterior derivative on differential forms
 extends, in a straightforward way,
 to an operator $\delta$ on $W$-valued
differential forms:

\begin{equation}
\delta a=dx^i\wedge \frac{\partial a}{\partial y^i}, \ \ \ \forall a\in \gm W\otimes \Lambda^* .
\end{equation}

By $\delta ^{-1}$, we denote its ``inverse" operator as defined  by:
\begin{equation}
\delta ^{-1}a=\frac{1}{p+q}y^{i} ( \frac{\partial }{\partial x^i}\per a)
\end{equation}
 when  $p+q>0$, and $\delta ^{-1}a=0$ when $p+q=0$.
Here    $a\in \gm W\otimes \Lambda^{q}$ is
  homogeneous of degree $p$ in $y$.

There is a  ``Hodge"- decomposition:
\begin{equation}
a=\delta \delta ^{-1}a+\delta ^{-1}\delta a+a_{0 0}, \ \ \ \forall a\in \gm W\otimes \Lambda^{*},
\end{equation}
where $a_{00}(x)$ is the  constant term of $a$, i.e, the $0$-form
term of $a|_{y=0}$ or  $a_{00}(x)=a(x, 0, 0, 0)$.
The operator $\delta$  resembles   most
basic properties of the usual exterior 
 derivatives.  For example,
$$\delta^{2}=0 \ \ \mbox{and } (\delta^{-1})^{2}=0 .$$

Let $\nabla$ be a torsion-free symplectic connection on $M$, and $\partial
: \gm W \lon \gm W \otimes \Lambda^{1}$ its induced covariant
derivative.

Consider   a connection on $W$ of the form:
\begin{equation}
\label{eq:connection}
D=-\delta +\partial + \ih [\gamma ,\cdot \ ],
\end{equation}
with $\gamma \in \gm W \otimes \Lambda^{1}$.

Clearly,  $D$ is a derivation with respect to the  Moyal-Weyl product, i.e.,
\begin{equation}
D(a* b)=a * Db +Da*  b.
\end{equation}

 A simple calculation yields that
\begin{equation}
D^{2}a=-[\ih \Omega , a] , \ \forall a\in \gm W,
\end{equation}
where
\begin{equation}
\label{eq:curvature}
\Omega =\omega -R+\delta \gamma -\partial \gamma -\ih \gamma^{2}.
\end{equation}
Here  $R=\frac{1}{4}R_{ijkl}y^{i}y^{j}dx^{k}\wedge dx^{l}$ and
$R_{ijkl}=\omega_{im}R^{m}_{jkl}$ is the curvature tensor
of the symplectic connection.

A connection  of the form (\ref{eq:connection}) 
 is called {\em Abelian} if
$\Omega $ is a scalar 2-form, i.e., $\Omega \in \Omega^{2}(M)[[\hbar ]]$.
It  is called  a {\em Fedosov
 connection } if it is Abelian and in addition $\gamma \in \gm W_{3}\otimes \Lambda^{1}$.
For an Abelian connection, the Bianchi identity  implies that
$d\Omega =D\Omega =0$, i.e., $\Omega \in Z^{2}(M)[[\hbar ]]$.
In this case,  $\Omega $ is called the Weyl curvature.

  \begin{thm} (Fedosov)
\label{label:f1}
Let  $\partial $ be   any torsion free symplectic connection,
 and
$\Omega =\omega +\hbar \omega_{1} +\cdots \in
Z^{2}(M)[[\hbar ]]$  a  perturbation of the symplectic form
in the space  $Z^{2}(M)[[\hbar ]]$.
 There exists
a unique $\gamma \in \gm W_{3}\otimes \Lambda^{1}$
such that $D$,  given by Equation (\ref{eq:connection}),
 is a Fedosov connection, which has Weyl curvature $\Omega$ and 
  satisfies
 $$\delta^{-1}\gamma =0.$$
\end{thm}
\pf It suffices to solve  the equation:
\begin{equation}
\omega -R+\delta \gamma -\partial \gamma -\ih \gamma^{2}=\Omega,
\end{equation}
which is equivalent to
\begin{equation}
\label{eq:r}
\delta \gamma =\tilde{\Omega}+\partial \gamma +\ih \gamma^{2}  ,
\end{equation}
where $\tilde{\Omega}=\Omega -\omega+R$.
Applying the operator $\delta^{-1}$ to Equation (\ref{eq:r})
and using the Hodge decomposition, we obtain
\begin{equation}
\label{eq:r-iteration}
\gamma =\delta^{-1}\tilde{\Omega}+\delta^{-1}(\partial \gamma +\ih \gamma^{2} ).
\end{equation}
Here we note that $\gamma_{00}=0$ since $\gamma$ is a $1$-form.

Since the operator
$\partial$ preserves the filtration and $\delta^{-1}$ raises it by
$1$, the iteration formula (\ref{eq:r-iteration}) has a
unique solution. Moreover
since $\delta^{-1}\tilde{\Omega}$ is at least  of degree 3, 
the solution $\gamma$ is indeed  in  $\gm W_{3}\otimes \Lambda^{1}$. \qed
{\bf Remark} The theorem above indicates that a Fedosov
connection is uniquely determined by a torsion free symplectic
connection $\nabla$ and a Weyl curvature $\Omega=\sum_{i}\hbar^{i}\omega_{i}\in
Z^{2}(M)[[\hbar ]]$. For this reason, we will say  that the connection $D$
defined above is a Fedosov connection corresponding to the pair
$(\nabla , \Omega )$.

If $D$ is a Fedosov connection,  the space of all parallel
sections $W_{D}$ automatically becomes  an  associative
algebra. Fedosov
proved that $W_{D}$ can be naturally
identified with $C^{\infty}(M)[[\hbar ]]$, and therefore
 induces a $*$-product on $C^{\infty}(M)[[\hbar ]]$, which we will
call  a {\em Fedosov $*$-product}.

Let $\sigma$ denote the projection from $W_{D}$ to
its center $C^{\infty}(M)[[\hbar ]]$ defined as
$\sigma (a)=a|_{y=0}$.

\begin{thm} (Fedosov)
\label{thm:f2}
$\sigma$ establishes an isomorphism between $W_{D}$ and
 $C^{\infty}(M)[[\hbar ]]$ as vector spaces. Therefore,
 it induces an associative algebra
structure on $C^{\infty}(M)[[\hbar ]]$, which is a $*$-product.
\end{thm}
\pf To prove that the two vector spaces are
isomorphic,
it suffices to show that for any $a_{0}(x, \hbar)\in C^{\infty}(M)[[\hbar ]]$ there is a unique
section $a\in W_{D}$ such that $\sigma (a)=a_{0}$.

The equation $Da=0$  can be written as

$$\delta a=\partial a +[\ih \gamma , a].$$
Applying the operator $\delta^{-1}$, it follows from the Hodge
decomposition that
\begin{equation}
\label{eq:parallel}
a=a_{0}+\delta^{-1}(\partial a+[\ih \gamma , a]).
\end{equation}

By the iteration method, we see that the equation above
has  a unique solution since $\delta^{-1}$ increases the filtration.
The rest of the claim can be easily verified. \qed

\section{Quantum exponential maps}

If $\nabla$ is flat and $\Omega=\omega$, the Fedosov connection
is  simply given by $D=-\delta +\partial $. In this case, the 
solution to Equation (\ref{eq:parallel}) can be expressed explicitly
as 
$$a=\sum_{k=0}^{\infty}\frac{1}{k!} (\partial_{i_{1}}\cdots
\partial_{i_{k}}a_{0})y^{i_{1}}\cdots y^{i_{k}}, $$
which is just the Taylor expansion of $exp_{x}^{*}a_{0}$ at
the origin. So the correspondence from $C^{\infty}(M)[[\hbar ]]$
to $W_{D}$ is indeed the pullback of the ($C^{\infty}$-jet  at
the origin of the) usual exponential map.
Thus for a general Fedosov connection, one may consider
the correspondence 
 $C^{\infty}(M)[[\hbar ]]\lon W_{D}$  as a 
``quantum correction" to the  exponential map.
In this section, we will  make this idea more
precise by introducing  the notion of quantum exponential
maps, which gives a simple characterization for  Fedosov connections.
As an application, we will realize geometrically
the equivalence between an arbitrary $*$-product and
one from Fedosov method, namely a Fedosov $*$-product.

\begin{defi}
A quantum exponential map is an $\hbar$-linear map
$\rho : \ \ C^{\infty}(M)[[\hbar ]] \lon  \gm W $ such that
\begin{enumerate}
\item $\rho (C^{\infty}(M)[[\hbar ]] )$ is a subalgebra of $\gm W$;
\item $\rho (a)|_{y=0}=a, \ \ \forall a\in C^{\infty}(M)[[\hbar ]]$;
\item $\rho (a)=a +\delta^{-1}da , \ \ \forall a\in C^{\infty}(M),  \ \ \mbox{mod } W_{2}$;
\item  $\rho (a)$ can be  expressed as a formal power series
in   $y$ and $\hbar$, with coefficients being derivatives of  $a$.
\end{enumerate}
\end{defi}

Given a quantum exponential map $\rho$, the Condition (ii)
implies that $\rho$ establishes an isomorphism
 between $C^{\infty}(M)[[\hbar ]]$
and its image as vector spaces.  Therefore,
$C^{\infty}(M)[[\hbar ]]$
becomes an associative algebra because of the first condition.
It is simple to see that the third condition  implies that
this is indeed a deformation of
the symplectic structure,  and the last one implies that it is
a $*$-product.

Clearly, for any Fedosov connection $D$, the
map from $C^{\infty}(M)[[\hbar ]]$ to $W_{D}$,
as constructed by Fedosov  (see Theorem \ref{thm:f2}),
satisfies all the properties
of a quantum exponential map. So quantum
exponential map always exists,  and one may 
consider  Fedosov construction as a
way of constructing a quantum exponential map.
In what follows, we will show that the converse
is also true. That is,

\begin{thm}
\label{thm:exponential}
Quantum exponential maps are equivalent to Fedosov connections.
\end{thm}

Before proving this theorem, we
start with investigating  the  following  closely related question:
what kind of subalgebras of $\gm W$ 
arises from  a  Fedosov connection? 

\begin{pro}
\label{pro:abelian}
Suppose that $\cala\subseteq  \gm W$ is   a subalgebra satisfying
the following  conditions:
\begin{enumerate}
\item ``completeness"- for  any $x_{0}\in M$, and  any $a(y, \hbar )\in W_{x_{0}}$,
there is an element $\tilde{a}(x, y, \hbar )\in \cala $ such that
$\tilde{a}(x_{0}, y, \hbar )=a(y, \hbar )$.
\item ``uniqueness"- if $\tilde{a}$ and $\tilde{b}
 \in \cala$ such that
$\tilde{a}|_{x_{0}}=\tilde{b}|_{x_{o}}$, then $\tilde{a}_{*}|_{x_{0}}
=\tilde{b}_{*}|_{x_{0}}$,
where $\tilde{a}$ and  $\tilde{b}$ are considered as maps: $M\lon W$,
 and $\tilde{a}_{*}$ and $\tilde{b}_{*}$ refer to their derivatives.
\end{enumerate}
Then, there exists a unique Abelian  connection $D$
such that $\cala  \subseteq W_{D}$.
\end{pro}
\pf Take a torsion free symplectic connection $\nabla$,
and let $\partial :\gm W\lon \gm W \otimes \Lambda^{1}$
be the corresponding covariant derivative.
For any $x_{0}\in M$, introduce an operator $\rho_{x_{o}}:
W_{x_{0}}\lon W_{x_{0}}\otimes \Lambda^{1}$ by
$$\rho_{x_{0}}(a(y, \hbar ))=(\delta -\partial ) \tilde{a}(x, y, \hbar )|_{x=x_{0}},$$
where $\tilde{a}(x, y, \hbar )\in \cala$ such that
$\tilde{a}(x_{0}, y, \hbar )=a(y, \hbar ) $.
By assumption, the map $\rho_{x_{0}}$  is  well-defined and  is  in
fact   a derivation of
the algebra $W_{x_{0}}$.
Therefore, there is a unique element $\gamma_{x_{0}}\in W_{x_{0}}\otimes
\Lambda^{1}$ with $\gamma_{x_{0}}|_{y=0}=0$ such that
$\rho_{x_{0}}=ad\ih \gamma_{x_{0}}=[\ih  \gamma_{x_{0}}, \cdot ]$.
Applying  this process pointwisely, we obtain  a global
section $\gamma \in W_{1}\otimes \Lambda^{1}$
with $\gamma_{0}\stackrel{def}{=}\gamma|_{y=0}=0$.
Let $D=-\delta +\partial +[\ih \gamma , \cdot ]$.
Then, $D$ is a connection on the Weyl bundle $W$ and
 satisfies the condition: $D \tilde{a}(x, y, \hbar )=0$
for all $\tilde{a}(x, y, \hbar )\in \cala$.
As in Section 2, let  $\Omega  =\omega -R+\delta \gamma -\partial \gamma -\ih \gamma^{2}$
 denote the Weyl
curvature of the connection $D$. It thus follows that
$$D^{2} \tilde{a}=[\ih \Omega , \tilde{a} ]=0,
  \ \ \ \forall \tilde{a}\in \cala. $$
Since $\cala$ is complete, $\Omega$
belongs to the center. Therefore, $\Omega \in Z^{2}(M)[[\hbar ]]$.
In other words, $D$ is an Abelian connection. \qed

The following lemma gives a simple sufficient condition for a subalgebra 
$\cala \subseteq \gm W$ being
complete.

\begin{lem}
\label{lem:complete}
Let $\cala \subseteq  \gm W$  be a subalgebra with unit.
Suppose that for any $a_{0}(x)\in C^{\infty}(M)$, there is
$\tilde{a}\in \cala $ such that
$$\tilde{a}=a_{0}+\delta^{-1}da_{0}, \ \ \ \ \mbox{mod } W_{2}. $$
Then, $\cala$ is complete.
\end{lem}
\pf We will  prove, by induction,
 the following  statement: $\forall x_{0}\in M  \  \mbox{ and }
a(y, \hbar )\in W_{k}|_{x_{0}}$, there is an
element  $\tilde{a}\in \cala$
such that  $\tilde{a}(x_{0}, y, \hbar )-a(y, \hbar ) \in  W_{k+1}|_{x_{0}} $.
The  conclusion will then follows immediately   from an 
iteration argument.

By assumption, the statement holds for both $k=0$ and $k=1$.

Assume that $a(y, \hbar )= \hbar^{j}y^{i_{1}}\cdots y^{i_{p}} $
with $2j+p=k$.
Now
\be
&&a(y, \hbar )-\hbar^{j}y^{i_{1}}*\cdots  *y^{i_{p}}\\
&=&\hbar^{j} (y^{i_{1}}\cdots y^{i_{p}}-y^{i_{1}}*\cdots  *y^{i_{p}}) \\
 &=& \sum_{2i +s=p} \hbar^{i+j}C_{i+j, j_{1} \cdots j_{s}}y^{j_{1}}\cdots
y^{j_{s}},
\ee
where $i\geq 1$, and  therefore $s=p-2i=k-2j-2i<k$.

Applying the induction assumption for each term $y^{j_{1}}\cdots
y^{j_{s}}$
in the above  summation,
we conclude that there is an element $\tilde{b}\in \cala$ such that
$$a(y, \hbar )-\hbar^{j}y^{i_{1}}*\cdots  *y^{i_{p}}=\tilde{b}(x_{0}, y, \hbar ), \ \ \ \mbox{mod } W_{k+1} .$$
On the other hand, for each $y^{i_{l}}$,  there
is $\tilde{a}^{i_{l}} \in \cala$ such that $\tilde{a}^{i_{l}}|_{x_{0}}=y^{i_{l}}$,
$\mbox{mod } W_{2}$, for $1\leq l\leq p$.
It is then clear that
$$ \hbar^{j}y^{i_{1}}*\cdots  *y^{i_{p}}=\hbar^{j}\tilde{a}^{i_{1}}*\cdots
*\tilde{a}^{i_{p}} |_{x_{0}}, \ \ \mbox{mod } W_{k+1}.$$
Thus, we have
$$a(y, \hbar )=\hbar^{j}\tilde{a}^{i_{1}}*\cdots
*\tilde{a}^{i_{p}} |_{x_{0}} +\tilde{b}(x_{0}, y, \hbar ),   \ \ \mbox{mod } W_{k+1}.$$
This concludes the proof.
\qed

Now we are ready to  formulate the following
 result, which  gives a  criterion  characterizing
a Fedosov algebra.

\begin{thm}
\label{thm:fedosov-algebra}
Let $\cala$ be a subalgebra of $\gm W$ satisfying:
\begin{enumerate}
\item $\forall a_{0}\in C^{\infty}(M)$, there is an element
$\tilde{a}\in \cala $ such that
$$\tilde{a}=a_{0}+\delta^{-1}da_{0}, \ \ \ \ \mbox{mod } W_{2} ;$$
\item  if $a$ and $b \in \cala$ such that
$a|_{x_{0}}=b|_{x_{o}}$, then $a_{*}|_{x_{0}}=b_{*}|_{x_{0}}$,
where $a$ and  $b$ are considered as  maps: $M\lon W$, and
$a_{*}$ and $b_{*}$ refer to their derivatives.
\end{enumerate}
Then, $\cala$ is a Fedosov algebra, i.e. $\cala$ arises
from a Fedosov connection.
\end{thm}
\pf According to Lemma \ref{lem:complete}, 
 $\cala$ is complete. Thus by Proposition
\ref{pro:abelian}, there is an Abelian connection
$D$ such that $\cala \subseteq  W_{D}$.

\begin{lem}
This  Abelian connection $D$ is in fact a Fedosov connection.
\end{lem}
\pf By assumption, for any $a_{0}\in C^{\infty}(M)$,
there is $\tilde{a}\in \cala $ such that
$$\tilde{a}=a_{0}+\delta^{-1}da_{0}, \ \ \ \ \mbox{mod } W_{2}. $$

It follows from $D\tilde{a}=0$ that
$$(\delta -\partial )\tilde{a}=[\ih \gamma , \tilde{a} ]. $$
The degree zero term of the LHS is easily seen to be zero,
while the degree zero term of the RHS is $\{\gamma_{1},  \delta^{-1}da_{0}\}$,
where $\gamma_{1}$ is the degree one term of $\gamma$.
It thus follows that  $\{\gamma_{1},  \delta^{-1}da_{0}\}=0$.
Since $a_{0}$ is arbitrary, $\gamma_{1}$ must be
constant with respect to $y$. However it is linear
in $y$, it has to be identically  zero.
This implies that $\gamma \in W_{2}\otimes \Lambda^{1}$.

Denote by $\gamma_{2}$ the degree 2 term of $\gamma$,  and assume
that $\gamma_{2}=\sum r_{ij, k}(x)y^{i}y^{j}dx^{k}$.
Note that this is the most general form  since
$\gamma_{0}=0$.
Also, we may always assume that $r_{ij, k}=r_{ji, k}$.
Since $D$ is Abelian, its  curvature $\Omega \in Z^{2} (M)[[\hbar ]] $.
Assume that
\begin{equation}
\label{eq:curvature2}
\Omega =\sum_{i=0}^{\infty}\hbar^{i}\omega_{i}
=\omega_{0}+\hbar \omega_{1}+\hbar^{2} \omega_{2}
+\cdots  .
\end{equation}

On the other hand, according to Equation (\ref{eq:curvature}), we have
\begin{equation}
\label{eq:curvature3}
\Omega = \omega -R+\delta \gamma -\partial \gamma -\ih \gamma^{2}.
\end{equation}

Comparing the  degree zero  terms of Equations
 (\ref{eq:curvature2}) and
(\ref{eq:curvature3}), it follows immediately
that $\omega_{0}=\omega $.
Now the terms $R$, $\partial \gamma $ and
$\ih \gamma^{2}$ are all of degree  not less than $2$, so the 
only degree $1$ term in Equation  (\ref{eq:curvature3})
would be  $\delta \gamma_{2}$.
Hence $\delta \gamma_{2}=0$. On the other hand,
a simple calculation yields that
\be
\delta \gamma_{2}&=&\sum ( r_{ij, k}(x)y^{j}dx^{i}\wedge dx^{k}
+r_{ij, k}(x)y^{i}dx^{j}\wedge dx^{k} )\\
&=&\sum  r_{ij, k}(x)y^{j}dx^{i}\wedge dx^{k}
+\sum r_{ji, k}(x)y^{j}dx^{i}\wedge dx^{k}\\
&=&2\sum  r_{ij, k}(x)y^{j}dx^{i}\wedge dx^{k} .
\ee
It thus follows that  for any  $i\neq k$,
$$ \sum_{j} r_{ij, k} (x) y^{j}=\sum_{j}r_{kj, i} (x) y^{j} . $$
That is, $ r_{ij, k}(x) = r_{kj, i} (x) $. Hence, $ r_{ij, k}$ is
completely symmetric with respect to $i, \ j, \ k$.

Let $\gm_{ijk}'=\gm_{ijk}+2r_{ij, k}$. Then, $\gm_{ijk}'$
defines  a torsion free symplectic  connection with induced
differential $\partial'=d+[\ih \gm ', \cdot ]$,
where $\gm '=\half \gm_{ijk}'y^{i}y^{j}dx^{k}$.
It is easy to see that   $\partial '=\partial +[\ih \gamma_{2}, \cdot ]$.
Let ${\gamma} '=\gamma -\gamma_{2}$.
Then ${\gamma} '\in W_{3}\otimes \Lambda^{1}$, and
$D=-\delta +\partial ' +[\ih \gamma ' , \cdot ]$. This
shows that $D$ is indeed a Fedosov connection.

\begin{lem}
For any $\tilde{a}\in W_{D}$, there is $\tilde{a}_{0}\in \cala$ and
$\tilde{b}\in W_{D}$ such that
$$\tilde{a}=\tilde{a}_{0}+\hbar \tilde{b}.  $$
\end{lem}
\pf  Let $a=\tilde{a}|_{y=0} =a_{0}(x)+\hbar a_{1}(x)+\cdots \in
C^{\infty}(M)[[\hbar ]]$.
Take $\tilde{a}_{0}\in \cala $ such that
$\tilde{a}_{0}=a_{0}+\delta^{-1} da_{0} \ \ \ \mbox{mod } W_{2}$,
which is always possible by assumption. Thus,
$\tilde{a}_{0}|_{y=0}=a_{0}(x)+O (\hbar )$ and then
$(\tilde{a}-\tilde{a}_{0})|_{y=0} =O(\hbar )$.
However, we know that  $\tilde{a}-\tilde{a}_{0} \in W_{D}$ since
 $\tilde{a}_{0} \in \cala
\subseteq  W_{D}$. This implies that $\tilde{a}-\tilde{a}_{0} =\hbar \tilde{b}$
for some  $\tilde{b}\in W_{D}$.
This concludes the proof of the lemma. \qed

By using the lemma above repeatedly, one  immediately
obtains the other inclusion, i.e.,
$W_{D}\subseteq \cala$.  This concludes our proof of theorem
\ref{thm:fedosov-algebra}.  \qed

Clearly, Theorem \ref{thm:exponential} is an immediate consequence of
Theorem \ref{thm:fedosov-algebra}. As  another application, we will  give a
geometric constructive proof for the following (see  \cite{BCG} \cite{Deligne}
\cite{NT1}): 

\begin{thm}
\label{thm:equivalent}
Any $*$-product is equivalent to  a Fedosov $*$-product.
\end{thm}
\pf As in Section 2, let  $P=TM$ be the     regular Poisson manifold
equipped with  the fiberwise linear symplectic structures.
Let $\calp=\gm (\cup_{m\in M}
\gm_{0}^{\infty}(T_{m}M, \reals ))$, where $\gm_{0}^{\infty}(T_{m}M, \reals )$
denotes the set of $\infty$-jets at $0$ of real valued functions
on $T_{m}M$.
Thus, $\calp$ is a Poisson algebra with a naturally induced
Poisson structure.

\begin{lem}
Any $*$-product on  $C^{\infty}(M)[[\hbar ]]$ induces a
$*$-product on the Poisson manifold $P=TM$ so  that
there is an algebra embedding $\rho :
C^{\infty}(M)[[\hbar ]]\lon \calp [[\hbar ]]$.
\end{lem}
\pf Take any torsion-free symplectic connection $\nabla$, and for
any fixed $m\in M$,
let $Exp_{m}$ be the formal symplectic exponential
map introduced by Emmrich and Weinstein \cite{EW}.
Then, $(Exp_{m})^{*}: \ C^{\infty}(M)\lon \gm_{0}^{\infty}(T_{m}M, \reals )$
is a Poisson algebra morphism, which in fact maps
$jet_{m}^{\infty}C^{\infty}(M)$ to $\gm_{0}^{\infty}(T_{m}M, \reals )$
isomorphically. Therefore,  any $*$-product on $C^{\infty}(M)[[\hbar ]]$
 induces a $*$-product
on $\gm_{0}^{\infty}(T_{m}M, \reals ) [[\hbar ]]$, hence on $C^{\infty}(T_{m }M) [[\hbar ]]$.
Thus, we obtain a \preferred  $*$-product (see Appendix B for the  definition),
 denoted as $\tilde{*}_{\hbar}$,
 on  the Poisson  manifold $P$,
and $Exp^{*}$ is clearly   an embedding of  the algebra. \qed

According to Proposition \ref{prop:main},  a  \preferred  $*$-product 
on $P$ is essentially unique.
 Hence, there exists an equivalence
operator: $T_{\hbar}=1+\hbar T_{1}+\hbar^{2}T_{2}+\cdots $, with
$T_{i}$ being leafwise  differential operators,  between
$(\calp [[\hbar ]], \tilde{*}_{\hbar})$  and  the
standard Weyl quantization $\gm W$.
Let $\cala =T_{\hbar} Exp^{*}(C^{\infty}(M)[[\hbar ]])\subset \gm W$.
Then $\cala $ is  a subalgebra of $\gm W$. It is simple
to see that $\cala$ satisfies all the conditions in Theorem
\ref{thm:fedosov-algebra}, so it is a Fedosov algebra.  This concludes 
the proof of the theorem. \qed



Since every  $*$-product is equivalent to  a Fedosov $*$-product,
its  {\em characteristic class}, as defined by Nest-Tsygan \cite{NT1},
  is  the class in $H^{2}(M)[[\hbar ]]$
of the Weyl curvature  of its equivalent Fedosov $*$-product.
This is well-defined since the Weyl curvatures of equivalent Fedosov $*$-products
are cohomologous (see \cite{F2}). In fact, Fedosov showed that two
Fedosov $*$-products are equivalent iff their Weyl curvatures are cohomologous \cite{F2}.
Thus (see \cite{NT1}),

\begin{thm}
\label{thm:equivalent1}
Two $*$-products are equivalent iff they have the same characteristic class.
\end{thm}

\section{Vey-$*$ products}

Let $(M, \omega )$ be a symplectic manifold, and  $\nabla $
 a torsion free symplectic  connection with covariant derivative
$\partial$. Define
$\partial  u=du$, $\partial^{2}u=\partial (\partial u)$,
$\partial^{k}u= \partial (\partial^{k-1} u)$ and so on.
It is simple to see that $\partial^{k}u$ is a
symmetric contravariant  $k$-tensor.
Let $\pi$ be the Poisson bivector field on $M$.
By $\pi^{k}$, we denote the $2k$-tensor:

$$(\pi^{k})^{i_{1}\cdots i_{k}, j_{1}\cdots j_{k}}=
\pi^{i_{1} j_{1}}\pi^{ i_{2} j_{2}}\cdots \pi^{ i_{k} j_{k}} . $$

Set
 $$P^{k}_{\nabla }(u, v)=<\pi^{k}, \partial^{k}u\otimes \partial^{k}v >
=\pi^{i_{1} j_{1}}\cdots \pi^{ i_{k} j_{k}}(\partial_{i_{1}}\cdots \partial_{i_{k}}u)
(\partial_{j_{1}}\cdots \partial_{j_{k}}v). $$
In particular, $P^{0}_{\nabla }(u, v)= uv$ and $P^{1}_{\nabla  }(u, v)
=\{u, v\}$.

\begin{defi} \cite{BFFLS} 
A Vey $*$-product is a star product on $C^{\infty}(M)[[\hbar ]]$
such that
\begin{equation}
u*_{\hbar}v=\sum \frac{1}{k!}(-\frac{i\hbar}{2})^{k} Q_{k}(u, v),
\end{equation}
where  $Q_{k}$ is a bidifferential operator
of maximum order $k$ in each argument
and its principal symbol coincides with that
of $P_{\nabla  }^{k}(u, v)$.
\end{defi}

The main result of the section is

\begin{thm}
\label{thm:vey}
Any Fedosov $*$-product is a Vey $*$-product.
Moreover, if  $D=-\delta +\partial +[\ih \gamma , \cdot ]$
as in Theorem \ref{label:f1}, then 
\begin{equation}
Q_{2}(u, v)=P^{2}_{\nabla}(u, v)+C_{2}(u, v),
\end{equation}
where $C_{2}$ is a bidifferential operator of maximum order $1$ in
each argument.
\end{thm}

Given a Fedosov connection $D=-\delta +\partial +[\ih \gamma , \cdot ]$,
any element
\begin{equation}
\label{eq:a}
a(x, y, \hbar )=\sum \hbar^{i}a_{i, j_{1} \cdots j_{k}}(x)
y^{j_{1}}\cdots y^{j_{k}}
\end{equation}
 in $W_{D}$
 is determined by the iteration
formula (\ref{eq:parallel}). Assume that $a_{0}(x)=\sigma (a)=a(x, 0, \hbar )\in C^{\infty}(M)$.
It is clear that each  coefficient in Equation (\ref{eq:a}) can be expressed
as  $a_{i, j_{1} \cdots j_{k}}(x)=\cald_{i, j_{1} \cdots j_{k}}a_{0} $
for certain differential operator $\cald_{i, j_{1} \cdots j_{k}}$.
We will say  that the term $\hbar^{i}a_{i, j_{1} \cdots j_{k}}(x)
y^{j_{1}}\cdots y^{j_{k}}$ is of  degree $(i+k, s)$ if
$\cald_{i, j_{1} \cdots j_{k}}$ is  a differential operator of
degree $s$.

\begin{pro}
\label{pro:order}
Under the same hypothesis as above, assume that
$$a(x, y, \hbar )=\sum \hbar^{i}a_{i, j_{1} \cdots j_{k}}(x)
y^{j_{1}}\cdots y^{j_{k}}  \in W_{D}. $$
Then
\begin{enumerate}
\item $a_{i, j_{1} \cdots j_{k}}(x)  =\cald_{i, j_{1} \cdots j_{k}}a_{0}$,
where $\cald_{i, j_{1} \cdots j_{k}} $ is a  differential operator
of degree not greater than $i+k$.
\item $$a(x, y, 0)=\sum_{k=0}^{\infty}\frac{1}{k!} (\partial_{i_{1}}\cdots
\partial_{i_{k}}a_{0})y^{i_{1}}\cdots y^{i_{k}}+
H, $$
where all terms in the remainder $H$   are of    degree $(l, k)$ with $l<k$.
\item  For any $j$, the order of the differential operator
  $D_{1, j} $ is not greater than $1$.
\end{enumerate} 
\end{pro}
\pf Since $a(x, y, \hbar ) $ is
generated by iteration formula (\ref{eq:parallel}), we prove (i) by induction.
For this purpose, we  need to analyses
the effect of the  operators $\delta^{-1}\partial $ and
$\delta^{-1}[\ih \gamma , \cdot ]$ on any term of degree
$(i+k, s)$. It is obvious that
$\delta^{-1}\partial $ maps terms of degree $(i+k, s)$
to those  of degree $(i+k+1, s+1)$.
On the other hand,
it is not difficult to check that
 $\delta^{-1}[\ih \gamma , \cdot ]$ maps terms of degree $(i+k, s)$
to those  of degree $(m, s)$ with $m\geq i+k$.
Claim (i) thus follows immediately.

To prove (ii), we need to analyses the $\hbar^{0}$-terms
produced from $\delta^{-1}[\ih \gamma , \cdot ]$
in  iteration formula (\ref{eq:parallel}).
The only possibility would be $ \delta^{-1}(-i \hbar \{\ih \gamma_{0}, \tilde{a}\})$, where
$\gamma_{0}=\gamma (x, y, 0, dx )$ and
$\tilde{a}=   a_{0, j_{1} \cdots j_{k}}(x)
y^{j_{1}}\cdots y^{j_{k}}$ is  assumed of degree $(l, k)$.
By Part (i), we know that  $l\leq k$.
Now $ \delta^{-1}(-i \hbar \{\ih \gamma_{0}, \tilde{a}\} )=
\delta^{-1}\{ \gamma_{0}, \tilde{a}\}$.
Since $\gamma_{0}$ is at least cubic in $y$ (deg $\gamma \geq 3$),
the latter is of degree $(l, k')$ with $k'$  being at least
$k+2$.
Therefore, all  terms of degree $(k, k)$ in $a(x, y, 0)$
come solely from the iteration of
the operator $\delta^{-1}\partial$ on elements of the same
form.
The conclusion thus follows immediately.

To prove (iii), we will concentrate on those terms of 
degree $(l, 2)$:
 \begin{equation} \label{eq:form1}
 \sum \hbar a_{1, j}(x)y^{j}.
\end{equation}
Since the operator $\delta^{-1}\partial $ maps any $(t, s)$-term
into that  of degree $(t+1, s+1)$, and
 $\delta^{-1}\partial $ does not produce any $\hbar$, the
only possible way to obtain such  a term is when
$\delta^{-1}\partial $ acts on those  of the form
$\hbar a_{1}(x)$. This however will never  happen since
we already  assume that $a(x, 0, \hbar )=a_{0}(x)$ is independent of $\hbar$.

As for the operator
  $\delta^{-1}[\ih \gamma , \cdot ]$, the only possible way  of generating
a  term of the form (\ref{eq:form1}) is
$\delta^{-1}(-i\hbar \{\ih \gamma' , a'\})=
\delta^{-1}\{\gamma ' , a'\}$,
where $\gamma'$ is the term in $\gamma$ having  the form:
$\hbar \gamma_{i, j}(x)y^{i}dx^{j}$ and $a'$ is the term in $a(x, y,  \hbar )$
of   the form:
$a'=a_{0, m}(x)y^{m} $.
Thus, $\delta^{-1}\{\gamma ' , a'\}=\hbar \gamma_{ij}(x) a_{0, m}(x)
\omega^{im}(x)y^{j} $ is of degree $(1, 2)$.
This concludes the proof. \qed

An  immediate consequence of Proposition \ref{pro:order}
is  the following:

\begin{cor}
Let $a(x, y, \hbar )\in W_{D}$, and let $x_{0}\in M$ be any point. Then,
$a(x_{0}, y, \hbar )=0$, $\forall y\in T_{x_{0}}M$ iff
$jet_{x_{0}}^{\infty} a(x, \hbar )=0$, where
$a(x, \hbar )=a(x, 0, \hbar )=\sigma (a( x, y, \hbar ) )$.
\end{cor}
\pf One  direction  follows directly from Part (1) of Proposition \ref{pro:order}.
To prove the other direction,
 let $a_{0}(x)=a(x, 0 ,0)$ and   assume that $a(x, \hbar )=
a_{0}(x)+\hbar a_{1}(x)+\hbar^{2}a_{2}(x) +\cdots $.
 We shall prove,
as the  first step, that $jet_{x_{0}}^{\infty} a_{0}(x)=0$.

Clearly, $a_{0}(x_{0})=0$. Assume that all derivatives
of $a_{0}$ up to the $k$-th order vanish at $x_{0}$.
According to  Part (ii) of Proposition \ref{pro:order},

$$a(x_{0}, y, 0)=
\sum_{k=0}^{\infty}\frac{1}{k!} (\partial_{i_{1}}\cdots
\partial_{i_{k}}a_{0})y^{i_{1}}\cdots y^{i_{k}}+
H,$$
where all terms in  $H$ are of degree $(l, k)$ with $l<k$.
Since the coefficient of $y^{i_{1}}\cdots y^{i_{k}}$
is  zero, it  follows that
$\frac{1}{k!} (\partial_{i_{1}}\cdots
\partial_{i_{k}}a_{0}) +(\cald a_{0} )(x_{0}) =0$, for
some differential operator $\cald$ of degree less than $k$.
By using the induction assumption, we deduce that
 $(\partial_{i_{1}}\cdots \partial_{i_{k}}) a_{0}(x_{0})=0$.
This proves that $jet_{x_{0}}^{\infty} a_{0}(x)=0$.
Let $a_{0}(x, y, \hbar)\in W_{D}$ be the parallel section
 corresponding  to
$a_{0}(x)$. Then $a_{0}(x_{0}, y, \hbar)=0, \ \forall y\in
T_{x_{0}}M$.
By considering the element  $\frac{1}{\hbar}[a(x, y, \hbar )- a_{0}(x_{0}, y, \hbar)]\in
W_{D}$, one yields that $jet^{\infty}_{x_{0}}a_{1}(x)=0$.
The conclusion thus follows  by using this argument repeatedly. \qed
%
{\bf Proof of Theorem \ref{thm:vey}}   Let $a_{0}(x) $ and $b_{0}(x)$ be
 any functions on $M$,  and $a(x , y, \hbar )$ and $b(x, y, \hbar )$
their  corresponding parallel sections  in $W_{D}$.
Assume that
$$a(x, y, \hbar )=\sum \hbar^{i}a_{i, j_{1} \cdots j_{k}}(x)
y^{j_{1}}\cdots y^{j_{k}},  \ \ \ \mbox{and} $$
$$b(x, y, \hbar )=\sum \hbar^{i}b_{i, j_{1} \cdots j_{k}}(x)
y^{j_{1}}\cdots y^{j_{k}}.  $$

By the  definition of Moyal-Weyl product,

\be
&&a(x, y, \hbar )*b(x, y, \hbar )|_{y=0}\\
&=&\sum (-\frac{i}{2})^{p} \frac{1}{p!} \hbar^{k+l+p}
a_{k, i_{1}\cdots i_{p}}(x)b_{l, j_{1}\cdots j_{p}} (x)\pi^{t_{1}s_{1}}\cdots
\pi^{t_{p}s_{p}}\frac{\partial (y^{i_{1}}\cdots y^{i_{p}} )}{\partial
y^{t_{1}} \cdots \partial y^{t_{p}} }\frac{\partial (y^{j_{1}}\cdots y^{j_{p}} )}{\partial
y^{s_{1}} \cdots \partial y^{s_{p}} } \\
&=&\sum_{n}\sum_{k+l+p=n}  (-\frac{i}{2})^{p} \frac{1}{p!} \hbar^{k+l+p}
a_{k, i_{1}\cdots i_{p}}(x)b_{l, j_{1}\cdots j_{p}}(x)\pi^{t_{1}s_{1}}\cdots
\pi^{t_{p}s_{p}}\frac{\partial (y^{i_{1}}\cdots y^{i_{p}} )}{\partial
y^{t_{1}} \cdots \partial y^{t_{p}} }\frac{\partial (y^{j_{1}}\cdots y^{j_{p}} )}{\partial
y^{s_{1}} \cdots \partial y^{s_{p}} } \\
&=&\sum_{n}\sum_{k+l+p=n}  (-\frac{i}{2})^{p} \frac{1}{p!} \hbar^{k+l+p}
(\cald_{k, i_{1}\cdots i_{p}} a_{0})( \cald_{l, j_{1}\cdots j_{p}} b_{0})
\pi^{t_{1}s_{1}}\cdots
\pi^{t_{p}s_{p}}\frac{\partial (y^{i_{1}}\cdots y^{i_{p}} )}{\partial
y^{t_{1}} \cdots \partial y^{t_{p}} }\frac{\partial (y^{j_{1}}\cdots y^{j_{p}} )}{\partial
y^{s_{1}} \cdots \partial y^{s_{p}} }
\ee
According to Proposition \ref{pro:order},
 the  order of the differential operator
$\cald_{k, i_{1}\cdots i_{p}}$ is not greater than $k+p$, which
is less than or equal to $n$, similarly for the order of
$\cald_{l, j_{1}\cdots j_{p}} $.
If $l\neq 0$, then $k+p=n-l<n$. The order of $\cald_{k, i_{1}\cdots i_{p}}$
is less than $n$. Similarly if $k\neq 0$, the order of
$\cald_{l, j_{1}\cdots j_{p}} $ is less than $n$.
So in order to have the maximum order in both arguments,
it is necessary that $k=l=0$, and in this case $p=n$.
Using Part (ii) of Proposition \ref{pro:order}, it
is simple to see that
the principal part is

$$ \sum (-\frac{i\hbar }{2})^{n} \frac{1}{n!} (\partial_{i_{1}}\cdots
\partial_{i_{n}}a_{0})(\partial_{j_{1}}\cdots \partial_{j_{n}}b_{0})
\pi^{i_{1}j_{1}}\cdots \pi^{i_{n}j_{n}} .$$

For the $\hbar^{2}$-term, we need  $k+l+p=2$.
If $k=l=0$ and $p=2$, we obtain the principal term.
When $k=p=1$ and $l=0$,  then
$a_{1, i_{1}}(x)=\cald_{1, i_{1}} a_{0}$
and $b_{0, j_{1}}(x)=\cald_{0, j_{1}}b_{0}$.
According to Proposition \ref{pro:order}, both $\cald_{1, i_{1}}$ and $\cald_{0, j_{1}}$
have degrees not greater than $1$;
similarly for the case that $l=p=1$ and $k=0$.

This concludes the proof. \qed

As an immediate consequence, we  obtain the following result, which was first
proved by Lichnerowicz using homological methods \cite{L}.

\begin{cor} (Lichnerowicz)
Any $*$-product  on a symplectic manifold is equivalent to a
Vey-$*$-product.
\end{cor}
{\bf Remark} As we know, the equivalence class of a
Fedosov $*$-product is determined by the class
of its Weyl curvature in $H^{2}(M)[[\hbar ]]$,
which is independent of the symplectic connection in the
construction. On the other hand, as we see in Theorem
\ref{thm:vey}, the symplectic connection is indeed completely
reflected in the Fedosov $*$-product itself. In fact
 if two Fedosov connections
$D_{i}=-\delta +\partial_{i}+[\ih \gamma_{i}, \cdot ]$
induce an identical $*$-product on $C^{\infty}(M)[[\hbar ]])$,
it is necessary that their symplectic connections
coincide, i.e., $\partial_{1}=\partial_{2}$.
We would like to conjecture that:

{\em Two Fedosov connections
$D_{i}=-\delta +\partial_{i}+[\ih \gamma_{i}, \cdot ]$ induce an
 identical $*$-product on $C^{\infty}(M)[[\hbar ]])$,
iff $D_{1}=D_{2}$.}

To prove this, one needs to prove that the Weyl curvatures
of $D_{1}$ and $D_{2}$ are not only cohomologous, but  indeed coincide,
or equivalently $W_{D_{1}}=W_{D_{2}}$. In other words. one
needs to decode the Weyl curvature (not only its class in $H^{2}(M)[\hbar ]$)
 from a Fedosov $*$-product directly.

We end this section  with  the following  inverse question:

\begin{quote}
{\bf Question} Is every Vey $*$-product necessary a Fedosov
$*$-product?
\end{quote}

\section{Quantization of lagrangian submanifolds}

Lagrangian submanifolds play a fundamental role in the study
of symplectic manifolds.  In fact, according to  the ``symplectic creed",
everything can be thought  of as a lagrangian submanifold \cite{Weinstein:1977}.
It is very natural to ask: what  should  lagrangian submanifolds
correspond  to under a deformation quantization?
This is what we aim to answer in this section.
\begin{lem}
\label{lem:lag}
Let $*_{\hbar}$ be the Moyal-Weyl product for  the symplectic
vector space $V\cong \reals^{2n}$ with the standard
symplectic structure. Suppose that  $L\subset V$ is a lagrangian
subspace. If both $f$ and $g\in C^{\infty}(M) $ vanish on $L$,
 so does $f*_{\hbar}g$.
\end{lem}
\pf Note that the Moyal-Weyl formula is independent of
the choice of linear coordinates.
Let us  choose a lagrangian subspace $K$ supplementary to
$L$, and linear coordinates $(q_{1}, \cdots , q_{n})\in L$,
$(p_{1}, \cdots , p_{n})\in K$ such that
$\omega (q_{i}, p_{j})=\delta_{ij}$.
Then the claim follows directly from  Equation  (\ref{eq:mw}). \qed

Let $L$ be a lagrangian submanifold. For any $x\in L$,
by $W_{x}^{L}$ we denote the subspace   of the Weyl algebra
$W_{x}$  consisting   of all  elements which  vanish
  when being  restricted to  $T_{x}L$.  According to Lemma \ref{lem:lag},
$W_{x}^{L}$ is a subalgebra.
Let $W^{L}=\cup_{x\in L}W_{x}^{L}$. Then $W^{L}$ is
a subbundle of the Weyl bundle. Similarly, we can
define $W^{L}_{1}$, $W^{L}_{2}$, etc. according to the
natural filtration in $W$.

We say an element $a\in W_{p}\otimes \Lambda^{q}$ is in the
subspace $(W_{p}\otimes \Lambda^{q})_{L}$ if for any
$v_{1}, \cdots , v_{q}\in T_{x}L$, $a(v_{1}, \cdots , v_{q})\in
W_{p}^{L}$. By $(W\otimes \Lambda )_{L}$, we  denote
the direct sum  $\oplus_{p, q}(W_{p}\otimes \Lambda^{q})_{L}$.
 It is clear that  both $\delta$ and $\delta^{-1}$
preserve the space $(W\otimes \Lambda )_{L}$.
More precisely,  $\delta $ maps $(W_{p}\otimes \Lambda^{q})_{L}$
into  $(W_{p-1}\otimes \Lambda^{q+1})_{L}$, while  $\delta^{-1}$
maps $(W_{p}\otimes \Lambda^{q})_{L}$
into  $(W_{p+1 }\otimes \Lambda^{q-1})_{L}$.

Let $\nabla $ be a  torsion-free  symplectic connection on $M$ such that $L$
is totally  geodesic. Then its induced differential
$\partial $ maps $(W\otimes \Lambda^{q})_{L}$ into
$(W\otimes \Lambda^{q+1})_{L}$.

\begin{pro}
\label{pro:lag}
As in Theorem \ref{label:f1}, if in addition, $L$ is  totally
geodesic with respect to the symplectic connection $\nabla$,
and is  lagrangian with respect to the Weyl curvature $\Omega$, i.e.,
$i^{*}\Omega =0$, where $i: L\lon M$ is the embedding,
then $\gamma $ belongs to $(W_{3}\otimes \Lambda^{1})_{L}$.
Therefore, the corresponding Fedosov connection
 $D$ preserves $(W\otimes \Lambda )_{L}$.
\end{pro}
\pf According to Equation (\ref{eq:r-iteration}), $\gamma $ is determined
by the iteration formula
\begin{equation}
\label{eq:r-iteration1}
\gamma =\delta^{-1} \tilde{\Omega }+\delta^{-1}(\partial \gamma +\ih \gamma^{2} ),
\end{equation}
where $\tilde{\Omega }=\Omega -\omega +R$.
It is easy to see that $R\in (W_{2}\otimes \Lambda^{2})_{L}$
since $L$ is a totally geodesic lagrangian
submanifold. Therefore, it follows that $\tilde{\Omega }\in
(W_{2}\otimes \Lambda^{2})_{L}$.
It is then clear, from the above iteration formula, that $\gamma \in
(W_{3}\otimes \Lambda^{1})_{L}$ since $(W\otimes \Lambda^{*})_{L}$
is preserved under all the operations involved in Equation (\ref{eq:r-iteration1}). \qed

Let $W_{D}^{L}$ denote the subspace of $W_{D}$ consisting of
 sections  whose restriction to $L$ belong to $W^{L}$.

As an immediate consequence, we have

\begin{cor}
Under the same hypothesis as in Proposition \ref{pro:lag},
$W_{D}^{L}$ is a subalgebra of $W_{D}$.
\end{cor}

By $C_{L}^{\infty}(M)$, we denote the space of smooth  functions on $M$ which
 vanish on $L$.

\begin{pro}
Under the same hypothesis as in Proposition \ref{pro:lag},
 $$  \sigma ( W^{L}_{D})=C^{\infty}_{L}(M)[[\hbar ]] ,$$
where $\sigma $
is the isomorphism between  $W_{D}$ and $C^{\infty}(M)[[\hbar]]$ as introduced in
Section 2.

Therefore, $C^{\infty}_{L}(M)[[\hbar ]]  $ is a subalgebra
of the Fedosov $*$-algebra $(C^{\infty}(M)[[\hbar ]], *_{\hbar})$.
\end{pro}
\pf According to Theorem \ref{thm:f2}, if $a\in W_{D}$ with
 $\sigma (a)=a_{0}$,
$a$ is determined by the iteration formula
$$a=a_{0}+\delta^{-1}(\partial a+[\ih \gamma , a]). $$
If $a_{0}\in C^{\infty}_{L}(M)[[\hbar]]$,  which means that
$a_{0}\in W^{L}$, it follows immediately that
$a\in W^{L}$. This is because   $\gamma \in (W_{3}\otimes \Lambda^{1})_{L}$,
and all the operators involved preserve the space  $(W\otimes \Lambda )_{L}$.

Conversely, if $a\in W_{D}^{L}$, it is obvious
that $a_{0}=\sigma (a)\in W^{L}$. \qed

\begin{ex} 
If $f: M\lon N$ is a symplectic diffeomorphism, its graph
$G_{f}=\{(x, f(x) )|x\in M\}$ is a lagrangian
submanifold of $S=M\times \bar{N}$. Moreover if $f^{*}\Omega_{N} =\Omega_{M}$,
 where $\Omega_{M}$  and $\Omega_{N}$ are  Weyl curvatures on $M$ and $N$,
respectively,
 $G_{f}$ is   lagrangian with respect to
$(\Omega_{M}, \overline{\Omega_{N}} )$.
Let $S$ be equipped with    a product symplectic
connection $\nabla \times \tilde{\nabla }$. It is easy
to see that $G_{f}$ is totally geodesic iff $\tilde{\nabla }=f_{*}\nabla $.
In this case, $C^{\infty}_{L}(S)$ is a subalgebra of
the corresponding Fedosov $*$-product $(C^{\infty}(S)[[\hbar ]], *_{\hbar})$.
This  implies  that $f_{*}$ is an algebra morphism
between   Fedosov $*$-algebras  $(C^{\infty}(M)[[\hbar]],  *_{\hbar})$
and  $(C^{\infty}(N)[[\hbar]],  *_{\hbar})$.
\end{ex}

The following well known lemma indicates that a totally
geodesic symplectic  connection always exists    for any given lagrangian
submanifold.

\begin{lem}
Given any lagrangian submanifold $L\subset M$, there always exists
a  torsion-free symplectic connection on $M$ such that $L$ is totally
geodesic.
\end{lem}
\pf First, take any torsion-free connection $\nabla $ on $M$
such that $L$ is totally geodesic.
Any other connection can be written
as
\begin{equation}
\tilde{\nabla}_{X}Y=\nabla_{X}Y +S(X, Y) , \ \ \ \forall X, Y\in \calx (M),
\end{equation}
where $S$ is a $(2, 1)$-tensor.
Clearly, $\tilde{\nabla}$ is torsion-free iff
$S$ is symmetric, i.e., $S(X, Y)=S(Y, X)$ for any $X, Y\in \calx (M)$.

$\tilde{\nabla}$ is symplectic iff $\tilde{\nabla}_{X}\omega =0$.
The latter is equivalent to
\begin{equation}
\omega (S(X, Y), Z)-\omega (S(X, Z), Y)=(\nabla_{X}\omega )(Y, Z).
\end{equation}

Let $S$ be the $(2, 1)$-tensor defined by
the equation:

\begin{equation}
\label{eq:s}
\omega (S(X, Y), Z)=\frac{1}{3}[(\nabla_{X}\omega )(Y, Z)+
(\nabla_{Y}\omega )(X, Z)] .
\end{equation}

Clearly, $S(X, Y)$,   defined  in this way,  is symmetric
with respect to $X$ and $Y$.  Now

\be
&&\omega (S(X, Y), Z)-\omega (S(X, Z), Y)\\
&=&\frac{1}{3}[(\nabla_{X}\omega )(Y, Z)+ (\nabla_{Y}\omega )(X, Z)]
-\frac{1}{3}[(\nabla_{X}\omega )(Z, Y)+(\nabla_{Z}\omega )(X, Y)]\\
&=&\frac{1}{3}[(\nabla_{X}\omega )(Y, Z)+(\nabla_{Y}\omega )(X, Z)
+(\nabla_{X}\omega )( Y, Z)+(\nabla_{Z}\omega )(Y, X)]\\
&=&(\nabla_{X}\omega ) (Y, Z),
\ee
where the last step follows from  the identity:

$$(\nabla_{X}\omega )(Y, Z)+(  \nabla_{Y}\omega )(Z, X)+
(\nabla_{Z}\omega )(X, Y)=0.$$

This means that $\tilde{\nabla}$ is a torsion-free symplectic connection.

>From  Equation (\ref{eq:s}),  it follows that
if $X$, $Y$ and $Z$ are all tangent
to $L$,  $\omega (S(X, Y), Z)=0$ since $L$ is totally geodesic with
respect to $\nabla$.
Hence, $S(X, Y)$ is tangent to $L$ whenever $X$, $Y$ are tangent to $L$.
In other words, $L$ is totally geodesic with respect to  $\tilde{\nabla}$. \qed

Since a 
 $*$-product is always equivalent to a
Fedosov $*$-product, as a consequence,
any lagrangian submanifold, under a deformation quantization,
 becomes a subalgebra after
some ``quantum correction". More precisely, we have

\begin{thm}
Let $*_{\hbar}$ be a $*$-product on a symplectic manifold
$(M, \omega )$ with characteristic class $[\Omega ]\in H^{2}(M)[[\hbar ]] $.
Suppose that $L$ is a lagrangian submanifold
such that $i^{*}[\Omega] \in H^{2}(L)[[\hbar ]]$ vanishes, where
$i: L\lon M$ is the embedding.
Then there exists  an operator: $T_{\hbar}=1+\hbar T_{1}+
\hbar^{2} T_{2}+\cdots $,  with $T_{i}$ being
differential operators  on $M$,  such that
$T_{\hbar}(C^{\infty}_{L}(M)[[\hbar ]]) $ is a
subalgebra of $(C^{\infty}(M)[[\hbar ]], *_{\hbar})$.
\end{thm}
\pf By assumption, $i^* \Omega $ is an exact two-form on $L$, i.e.,
$i^* \Omega =d\theta_{L}$ for some $\theta_{L}\in \Omega^{1}(L)$.
Extending $\theta_{L}$ to a one-form on $M$, we may
assume that $i^* \Omega =d i^{*} \theta$ for some $\theta\in
\Omega^{1}(M)$. Let $\tilde{\Omega}=\Omega -d\theta$, and
 take a torsion-free symplectic connection $\nabla $ such that
$L$ is totally geodesic.
Let $\bar{*}_{\hbar}$ be the corresponding
Fedosov $*$-product  with  Weyl curvature $\tilde{ \Omega }$.
Then $ C^{\infty}_{L}(M)[[\hbar ]] $ is a  $\bar{*}_{\hbar}$-subalgebra.
According to Theorem \ref{thm:equivalent1}, $*_{\hbar}$ and  $\bar{*}_{\hbar}$  are equivalent
$*$-products. The conclusion thus follows immediately. \qed
{\bf Remark}\ \  (1). The quantum counterparts of lagrangian submanifolds,
according to Lu \cite{Lu}, are left ideals. However, it is not
clear how this can be realized for  $*$-algebras in our case.
It seems that a possible candidate  would be the space 
 $ C^{\infty}_{L}(M)[[\hbar ]] $ 
modified in a certain way. 

(2).  For a   symplectic manifold $M$, a coisotropic submanifold  
  is a submanifold  $C$
 such that the space of 
 functions vanishing  on $C$ becomes a Poisson
subalgebra. It is natural to expect  that the above result
can be generalized to any coisotropic submanifolds. But we cannot
prove this  at the moment because it is not clear if there
always exists a symplectic connection such that $C$ is
totally geodesic for a general coisotropic
submanifold $C$. 

\section{Quantum momentum maps}

This section is devoted to the study of deformation quantization of
a symplectic $G$-space. In particular, we will introduce the notion
of  quantum momentum maps, which plays the role of a quantum
analogue of the usual momentum maps.

Let $(M, \omega )$ be a symplectic $G$-space with   action
 $\Phi_{g}: M\lon M$, $\forall g\in G$.   A $*$-product on $M$ is called
{\em $G$-\equivariant} if for any $u$ and $v\in C^{\infty}(M)[[\hbar]]$,
\begin{equation}
\Phi^{*}_{g}(u*_{\hbar}v)=(\Phi^{*}_{g}u)*_{\hbar}(\Phi^{*}_{g}v).
\end{equation}

In general,  $M$  does not necessary   admit a $G$-\equivariant $*$-product. 
It is known \cite{BFFLS} \cite{L}  that the 
existence of such a $*$-product is closely related to the
existence of a $G$-invariant connection on the
manifold. More precisely,

\begin{pro}
\label{pro:natural} 
Let $M$ be a symplectic $G$-space. $M$ admits a $G$-\equivariant 
\nat  $*$-product
iff there exists a $G$-invariant connection on $M$,
where by a \nat $*$-product, we mean a $*$-product:
$$u*_{\hbar }v=\sum_{k} \hbar^{k} C_{k}(u, v),$$
where the $\hbar^2$-term $C_{2}$ is a bidifferential operator of
order $2$ in each argument.
\end{pro}
\pf  Assume that there exists a $G$-\equivariant \nat     $*$-product:
$$ u*_{\hbar }v=uv-\frac{i\hbar }{2}\{u , v\} + \half (\frac{i\hbar }{2})^{2}Q_{2}(u, v)+\cdots . $$
According to Proposition \ref{pro:nat-split}, there is a unique symplectic connection $\nabla$
such that
$$Q_{2}(u, v)=P_{\nabla}^{2}(u, v)+H(u, v), $$
where $H(u, v)$ is a bidifferential operator of maximum order  $1$ in each argument.
It is then clear that $\nabla $ is   $G$-invariant.

Conversely, suppose that there exists  a $G$-invariant connection
on $M$.  Using  the standard method, we can  always make it into
a $G$-invariant torsion free  symplectic connection $\nabla $.
Then, the corresponding Fedosov $*$-product (with  symplectic connection
$\nabla $, and  Weyl curvature $\omega $) will be $G$-equivariant. \qed

A $G$-invariant connection always exists if $G$ is compact. However, when
$G$ is non-compact, there exist some cases where $G$-invariant connections do not
exist. Various attempts have been made  in order to deal with such a situation.
For details, readers can consult  \cite{AC1} \cite{AC2} \cite{FS}.

In what follows,  nevertheless, we will  always assume that $ *_{\hbar}$
is a $G$-\equivariant $*$-product.
Then, the corresponding infinitesimal  action
 $\xi \lon \hat{\xi}$ defines  a
Lie algebra homomorphism from  $\frakg$ to the  derivation
space $Der C^{\infty}(M) [[\hbar ]]$ of the $*$-algebra.

By $\pol [[\hbar ]]$, we denote the space of formal power series of
$\hbar$ with coefficients in the universal
enveloping algebra $\pol$. Let $\frakg_{\hbar}$ be the deformed
Lie algebra:  $\frakg_{\hbar}= \frakg [[\hbar]]$ with the bracket
being defined as
\begin{equation}
[X, Y]_{\hbar}=-i\hbar [X, Y], \ \ \ \forall X, Y\in \frakg [[\hbar]].
\end{equation}
Then $\pol [[\hbar ]]$ can be identified with the universal enveloping algebra
of $\frakg_{\hbar}$, and therefore inherits an associative algebra structure.

\begin{defi}
A quantum momentum map is a homomorphism of associative
algebras:
$$\mu_{\hbar}: \pol [[\hbar ]] \lon C^{\infty}(M) [[\hbar ]], $$
such that for any $\xi \in \frakg$,
\begin{equation}
\label{eq:moment}
\hat{\xi} =\ad {\ih \mu_{\hbar}(\xi )},
\end{equation}
where both sides are considered as derivations on $C^{\infty}(M) [[\hbar ]]$.
\end{defi}

It is obvious, from definition, that
a necessary condition for the existence of  a momentum map
is that the derivation $\rho_{\xi}f=\hat{\xi}f$,
$\forall \xi \in \frakg$   be  inner.
Let us first assume    that this is true.
Thus there is a linear map  from $\frakg$ to
$C^{\infty}(M)[[\hbar ]]$, denoted by
$\xi \lon a_{\xi}$,  such that
for any $f\in C^{\infty}(M)$ and $\xi \in \frakg$,
\begin{equation}
\label{eq:moment1}
\hat{\xi}f=[\ih a_{\xi }, f].
\end{equation}

Therefore,
\be
\widehat{[\xi , \eta ]}f &=&[\hat{\xi }, \hat{\eta }]f\\
&=& [\ih  a_{\xi }, [\ih a_{\eta }, f]]-[\ih a_{\eta  }, [\ih a_{\xi }, f]]\\
&=&[[\ih a_{\xi }, \ih a_{\eta }], f].
\ee
On  the other hand, by definition,
$$\widehat{[\xi , \eta ]} f =[\ih a_{[\xi ,\eta ]}, f]. $$

Therefore,
$$[a_{[\xi , \eta ]}-\ih [a_{\xi }, a_{\eta } ], f]= 0,  \ \ \ \forall f\in
 C^{\infty}(M). $$

So   $a_{[\xi , \eta ]}-\ih [a_{\xi }, a_{\eta } ]$, as a function in
 $C^{\infty}( M)[[\hbar ]]$, 
is  constant.

Define $\lambda : \wedge^{2}\frakg \lon \complex [[h]]$ by
\begin{equation}
\label{eq:2cocycle}
\lambda (\xi ,\eta )=a_{[\xi , \eta ]}-\ih [a_{\xi }, a_{\eta } ], \ \ \ \forall \xi , \eta \in \frakg .
\end{equation}

\begin{pro}
\begin{enumerate}
\item $\lambda$ is a  Lie algebra 2-cocycle;
\item its   cohomology class $[\lambda ]\in H^{2}(\frakg , \complex [[\hbar]])
\cong H^{2}(\frakg )\otimes \complex [[\hbar ]] $ is independent
of the choice of the linear map $a_{\xi}$;
\item quantum momentum map exists iff $[\lambda ]=0$.
\end{enumerate}
\end{pro}
\pf  Assertions (i)-(ii) are quite obvious,  and left for  the  reader
to check.

For (iii), suppose that a  quantum momentum map $\mu_{\hbar}$
exists.  Then we may take $a_{\xi}=\mu_{\hbar}({\xi})$
as our linear map.  In this case,

\be
\lambda (\xi, \eta )&=&a_{[\xi , \eta ]}-\ih [a_{\xi}, a_{\eta }]\\
&=&\mu_{\ih}{[\xi , \eta ]} -\ih [ \muh \xi , \muh {\eta}]\\
&=&\mu_{\ih}([\xi , \eta ]-\ih  [{\xi} ,{\eta}]_{\hbar} ) \\
&=&0.
\ee

Conversely, if $[\lambda ]=0$, by adding a suitable coboundary we can always choose a
linear map $a:\frakg \lon C^{\infty}(M)[[\hbar ]]$
such that  Equation (\ref{eq:moment1}) holds and
$a_{[\xi , \eta ]}-\ih [a_{\xi }, a_{\eta } ]=0$.
In other words, $a$ is a Lie algebra
homomorphism from $\frakg_{\hbar}$
to the commutator Lie algebra of $C^{\infty}(M)[[\hbar ]]$.
Therefore, it  extends to an associative algebra
morphism:
$$\muh : \ U\frakg_{\hbar}\lon (C^{\infty}(M)[[\hbar ]] ,
*_{\hbar}). $$

This concludes the proof. \qed

According to Theorem \ref{thm:derivation}, derivations are automatically inner if $H^{1}(M)=0$.
Thus we have

\begin{thm}
\label{eq:mom-exist}
There exists a quantum momentum
map if  $H^{1}(M)=0$ and $H^{2}(\frakg )=0$. In particular,
   a quantum momentum map exists if $M$ is simply connected
 and $\frakg$ is semi-simple.
\end{thm} 

We note that the above condition is exactly the same sufficient
condition for the existence of a  classical momentum map  \cite{AM}.
However, there are many cases where classical momentum maps still exist
even if  this condition is no  longer satisfied.
It is reasonable to    expect that this phenomenon  would happen
 for quantum momentum maps as well. However,
we do not know too many examples except for the following:

\begin{ex}
Suppose that  $Q$ is a $G$-manifold with action $\phi_{g}$,
which admits  a $G$-invariant torsion free connection
$\nabla$. Let $ M=T^{*}Q$ be  equipped with the standard cotangent
bundle  symplectic structure.
The $G$-action naturally lifts to  a symplectic action  
$\Phi_{g}$ on $ M=T^{*}Q$ 
with  an equivariant  momentum map  $J : T^{*}Q \lon \frakg^{*}$\cite{AM}:
$$ <J (\xi_{q}), X> =<\hat{X} , \xi_{q}>, \ \  \ \xi_{q}\in T_{q}Q , \forall X\in \frakg ,$$
where $\hat{X}$ is the vector field on $Q$ generated by $X\in \frakg $.

To any differential operator $D$, we assign its (complete)   symbol
 as the polynomial $S_{D}$
on $T^{*}Q$ given by
\begin{equation}
\label{eq:symbol}
S_{D} (\xi_{q} )=D e^{<exp_{q}^{-1}x , \xi_{q}>}|_{x=q}, \ \ \ \ \forall \xi_{q}\in  T_{q}Q,
\end{equation}
where $exp_{q}: T_{q}Q\lon Q$ is the usual exponential map, defined in  a
neighborhood of $0$, corresponding to the connection $\nabla$.

This   assignment in fact establishes  an isomorphism
between the space $\cald$  of differential operators and  that of polynomials on $T^{*}Q$.

Deform the Lie bracket structure on $\calx (Q)$ according to:
\be
&& Z_{\hbar}f= -i \hbar Zf, \ \forall Z\in \calx (Q) \mbox{ and } f\in C^{\infty}(Q); \mbox{ and } \\
&&[Y, Z ]_{\hbar} =-i \hbar [Y, Z], \forall Y, Z\in \calx (Q).
\ee
Since differential operators  are generated by $C^{\infty}(Q)$ and 
$\calx (Q)$
over the module $C^{\infty}(Q)$, this deformed bracket induces
  an $\hbar$-depending multiplication on $\cald $\footnote{A
 more intrinsic viewpoint
is to think the tangent bundle $TQ$ as  a Lie algebroid, and the above 
construction as
 deforming  the algebroid structure by multiplying the factor $-i \hbar$. Then
such a  construction admits an immediate generalization to algebroids, which
should give rise to a $*$-product for the  Lie-Poisson structure associated
 with a  Lie algebroid.}, which  in turn induces a $\hbar$-depending
multiplication on the space of polynomials on $T^{*}Q$, hence a $*$-product
on $T^{*}Q$.
It is simple to see that for any $D\in \cald$,
$$ \Phi_{g}^{*} S_{D} =S_{g^{-1} \cdot D}, $$
where $g\cdot D\stackrel{def}{=}\phi_{g^{-1}}^{*}\smalcirc D \smalcirc \phi_{g}^{*}$.
In the equation above, by letting $g=\exp{tX}$,   $ \forall X\in \frakg$, and
taking derivative at $t=0$, one obtains immediately that
$$\hat{X}(S_{D})=S_{[\hat{X}, D]}, $$
where $\hat{X}$ on the LHS refers to the vector field on $Q$, while $\hat{X}$ on
the RHS stands  for the one on $T^{*}Q$ corresponding to $X\in \frakg$.
This  implies that $$\hat{X} f=[\ih J^{*}l_{X}, f], \ \ \ \forall f\in C^{\infty}(T^{*}Q).$$
In other words,
 $\muh X=J^{*} l_{X}$ defines  a quantum momentum map. However, it is not clear
how to express  $\muh f$  explicitly for a  general $f\in \pol [[\hbar]]$.
\end{ex}

When quantum momentum maps exist, they are, in general, not unique, as
in the classical  case. 

Assume that both  $\mu_{\hbar}$ and $\nu_{\hbar}: \pol [[\hbar ]] \lon C^{\infty}(M) [[\hbar ]]$
are quantum momentum maps. Let
$\tau_{\hbar}: \frakg \lon C^{\infty}(M) [[\hbar ]] $ be the map defined as:
$$\tau_{\hbar}(\xi )= \mu_{\hbar} (\xi )-\nu_{\hbar} (\xi  ), \ \ \ \ \forall \xi\in \frakg .$$
Then, for any $f \in C^{\infty}(M)$, $\ih [ \tau_{\hbar}(\xi ) , f]=0$. Thus it follows that
$\tau_{\hbar}(\xi )$ is a constant formal polynomial of
 $\hbar$. Also, it is easy to see
that
$$\tau_{\hbar} ([\xi , \eta ])=0. $$
That is,   $\tau_{\hbar}: \frakg \lon \complex    [[\hbar ]]$
is a 1-cocycle. 
 Since all 1-coboundaries are trivial, it follows that the
 quantum momentum map is unique  if $H^{1}(\frakg )=0$.

\begin{pro}
If $H^{1}(\frakg )=0$, then the quantum momentum map is unique.
\end{pro}

 In general, for $\xi\in \frakg$, we have $\mu_{\hbar}\xi = \nu_{\hbar} \xi +\tau_{\hbar} \xi$.
However, it is not clear how to express $\mu_{\hbar}f$,
 for  a general $f\in \pol [[\hbar ]]$,
in terms of $\nu_{\hbar}$ and $\tau_{\hbar}$.

To see the relation between  a  quantum momentum map  and
a  classical one, we start with
the following simple

\begin{lem}
\label{lem:qcm} Suppose that $$f*_{\hbar}g=\sum_{k}\hbar^{k} C_{k}(f, g)$$ is a $*$-product on a
symplectic manifold $M$.  Let $X\in \calx (M)$ be
a vector field on $M$, which is an inner derivation
when being considered as an operator on  $C^{\infty}(M)[[\hbar ]]$.
Assume that $X f=[\ih a , f]$,   $\forall f\in C^{\infty}(M)[[\hbar ]]$.
Then, modulo a constant formal  polynomial of $\hbar$,  $a$ is of the form $a=\sum_{i}\hbar^{2i}a_{2i}$,
where $a_{0}$ is a hamiltonian function generating the
vector field $X$,  and $a_{2k}$, for $k\geq 1$, is determined by
the equation:
$$ \{a_{2k}, f\}= \sum_{i+j=k-1, i\geq 0, j\geq 0} -2 \ih C_{2j+3}(a_{2i}, f).
$$
\end{lem}
\pf It is a direct verification, and  is left for the reader. \qed


As a vector space,  $\pol [[\hbar ]]$ is canonically isomorphic to
$\poll [[\hbar ]]$, the space of formal power series of $\hbar$ with
coefficients being polynomials on $\frakg^*$.
The isomorphism is established  by the symmetrization (see \cite{B} for
details). Therefore, the algebra structure on $\pol [[\hbar ]]$
induces a $*$-product on $\poll [[\hbar ]]$, which gives rise
to a deformation quantization for the Lie-Poisson structure $\frakg^*$. Below, we will identify
these two spaces and use them interchangeably if there is no confusion.

\begin{pro}
\label{pro:equ-moment}
Suppose that $*_{\hbar} $ is
a $G$-equivariant $*$-product on a symplectic manifold $M$. Assume that $\muh :  \pol [[\hbar]]
 \lon C^{\infty}(M)[[\hbar ]]$ is  a quantum momentum
map.
Then, $M$ is a hamiltonian $G$-space, i.e.,  the symplectic
$G$-action admits
 an equivariant  (classical) momentum  map $J$.
Moreover,
$$\muh f=J^{*}f +O(\hbar ),  \ \ \ \ \forall f\in \poll .$$
\end{pro}
\pf  Since $\muh \xi, \ \forall \xi\in \frakg$,  depends on $\frakg$
linearly, it defines a map
 $J:X \lon \frakg^*$  uniquely by
 the  relation:
$$\mu_{\hbar}(l_{\xi })=J^{*}l_{\xi }+ O(\hbar) , \ \ \forall \xi\in \frakg .$$
Then clearly $J$ is a (classical) momentum  map  according to
Lemma  \ref{lem:qcm}.

Write the  $2$-cocycle $\lambda$   defined by  Equation (\ref{eq:2cocycle}) as
$$\lambda =\lambda_{0}+\hbar \lambda_{1}+\cdots .  $$
Then each $\lambda_{i}: \wedge^{2}\frakg \lon \complex$ is a 
Lie algebra $2$-cocycle. Moreover,
it is simple  to see that
$$\lambda_{0}(\xi, \eta )= J^{*}l_{[\xi ,\eta ]}-\{J^{*}l_{\xi} , \ J^{*}l_{\eta }\} .$$
The vanishing of $\sigma$ implies that $\sigma_{0}=0$, which
means that $J$ is equivariant. The rest of the proposition
thus follows trivially. \qed

It is,  however,  not clear whether the converse 
of Proposition \ref{pro:equ-moment} is true or  not.
 We end this section by posing the following

\begin{quote}
 {\bf QUESTION} Does the existence of a classical moment
map imply the existence of a quantum moment map?
\end{quote}

\section{Quantum dual pair}

This is  a continuation of the last section. We will assume the same
hypothesis as in the previous  section, and in particular, assume  that
a quantum momentum  map  $\mu_{\hbar}: \pol [[\hbar ]]  \lon C^{\infty}(M) [[\hbar ]]$
exists.

Let $C^{\infty}(M)^{G}$ denote the space of all $G$-invariant
functions on $M$.

\begin{pro}
\label{prop:comm1}
$$(\muh \pol [[\hbar ]] )' \cong C^{\infty}(M)^{G}[[\hbar ]].$$
\end{pro}
\pf Let $f=\sum_{i}\hbar^{i} f_{i} \in  C^{\infty}(M)[[\hbar ]]$.
 Then $f$ commutes with  $\muh \pol [[\hbar ]]$ iff it commutes with
its  generators. That is,
 $[ f , \muh \xi ]=0, \ \ \forall \xi \in \frakg$.
The latter is  equivalent to that $\hat{\xi}f=0$, or in
other words, $f$ is $G$-invariant. \qed

According to  this proposition,  we have $\muh (\pol [[\hbar ]])\subseteq
(C^{\infty}(M)^{G}[[\hbar ]])'$.  In order to  describe the \commutant
$(C^{\infty}(M)^{G}[[\hbar ]])'$ completely, we need to extend the
quantum momentum  map $\muh$.

It is well-known that the star-product on $\poll [[\hbar ]]$  naturally extends
to a star-product on  smooth
functions $C^{\infty}(\frakg^{*})[[\hbar ]]$.
Below we will indicate that a quantum momentum  map, if it exists,
 extends to  $C^{\infty}(\frakg^{*})[[\hbar ]]$ as well. More
precisely, we have

\begin{pro}
\label{prop:extension}
Let  $\mu_{\hbar}: \poll [[\hbar ]]  \lon C^{\infty}(M) [[\hbar ]]$ be a quantum
momentum  map. Then it  naturally  extends  to an algebra morphism,
denoted by the same notation $\mu_{\hbar}$, from $C^{\infty}(\frakg^{*})[[\hbar ]]$
to $ C^{\infty}(M) [[\hbar ]]$ such that
\begin{enumerate}
\item for any $f\in C^{\infty}(\frakg^{*})$,  $\mu_{\hbar}f$
commutes with $C^{\infty}(M)^{G}[[\hbar ]]$, and
\item for any $f\in C^{\infty}(\frakg^{*})$,
$\muh f=J^{*}f +O(\hbar )$.
\end{enumerate}
\end{pro}
\pf  Given a smooth function $f\in C^{\infty}(\frakg^{*})$,
 for any $x_{0}\in M$, let $\tilde{f}_{x_{0}} (u), \ u\in \frakg^{*} $
 denote its  Talyor expansion
 at the point $u_{0}=J(x_{0} )$. Define\footnote{We are grateful to
A. Weinstein for suggesting this method of extension.}

\begin{equation}
(\mu_{\hbar} f )(x_{0} )=(\mu_{\hbar} \tilde{f}_{x_{0}} )|_{x=x_{0}}.
\end{equation}

It is clear that  this definition  coincides with the   original
 $\mu_{\hbar}$ when $f$ is a polynomial, so it is indeed
an extension of the given quantum momentum map. It is
also obvious, from definition, that  $\mu_{\hbar}f$
commutes with $C^{\infty}(M)^{G}[[\hbar ]]$.

It is simple to see that each term in the
expansion of  $\tilde{f}_{x_{0}}-f(J(x_{0}))$,
 as a function in $u$,
 is a homogeneous polynomial in $u-u_{0}$.
Hence, we have  $\muh f=J^{*}f +O(\hbar )$.

It remains to check that $\mu_{\hbar}$ is an algebra
homomorphism. This follows from the fact that
$\widetilde{(f*_{\hbar }g)}_{x_{0}}=
\tilde{f}_{x_{0}} *_{\hbar }\tilde{g}_{x_{0}}$,
for all $f$ and $g\in C^{\infty}(\frakg^*  )$.

To see this, we write $f*_{\hbar }g=\sum \hbar^{k}C_{k}(f, g)$.
Then, this essentially  follows from
the fact that $\widetilde{C_{k}(f, g)}_{x_{0}}=C_{k}(\tilde{f}_{x_{0}} , \tilde{g}_{x_{0}})$, which is in turn a consequence of the fact that
$\poll  [[\hbar ]]$ is a subalgebra under the $*$-product on $C^{\infty}(\frakg^{*})[[\hbar ]]$. \qed
\begin{pro}
\label{prop:comm3}
Under the same hypothesis as in Proposition \ref{prop:extension}, if,
in addition,  the $G$-action is free and proper,
$$(C^{\infty}(M)^{G}[[\hbar ]])'\cong \muh
(C^{\infty}(\frakg^{*})[[\hbar ]] ). $$
\end{pro}
\pf  According to Proposition \ref{prop:extension}, $ \muh  (C^{\infty}(\frakg^{*})[[\hbar ]] )
\subseteq (C^{\infty}(M)^{G}[[\hbar ]])'$.
To prove the other direction,  let us assume that
$f=\sum_{i}\hbar^{i} f_{i} \in (C^{\infty}(M)^{G}[[\hbar ]])'$.

It is clear that for any $g\in C^{\infty}(M)$,
$$\ih [f , g] =\{f_{0}, \ g\}+O(\hbar ). $$
Then  it follows that
$\{f_{0}, \ g\}=0$,
 $ \forall \ g\in C^{\infty}(M)^{G}$.
 This  implies that
$f_{0}= J^{*}f'_{0} $ for some smooth function $f'_{0}\in C^{\infty}(\frakg^{*})$.  Now
according to  Proposition \ref{prop:extension},
$\muh f'_{0}=J^{*}f'_{0} +O(\hbar )
=f_{0}+O(\hbar )$. Therefore, $f-\muh f'_{0} =\hbar \tilde{f}$,
where $\tilde{f}\in  (C^{\infty}(M)^{G}[[\hbar ]])' $ since
both $f$ and $\muh f'_{0}$ belong to $ (C^{\infty}(M)^{G}[[\hbar ]])'$.
By repeating  the same argument on $\tilde{f}$ and  so on, we  deduce that
$f\in \muh (C^{\infty}(\frakg^{*})[[\hbar ]] )$. \qed

Combining Propositions \ref{prop:comm1}-\ref{prop:comm3},  we have

\begin{thm}
\label{thm:dual-pair}
Suppose that $*_{\hbar}$ is a  $G$ equivariant $*$-product on $C^{\infty}(M)[[\hbar ]]$
with a quantum momentum map $\muh$.
Assume that the action is free and proper, then
\be
(C^{\infty}(M)^{G}[[\hbar ]])' &\cong & \muh
(C^{\infty}(\frakg^{*})[[\hbar ]] )\\
(  \muh C^{\infty}(\frakg^{*})[[\hbar ]] )'&\cong
&C^{\infty}(M)^{G}[[\hbar ]] .
\ee
\end{thm}

Recall that if $A$ is an associative algebra and $B\subseteq A$ is a 
subset, then $B'$, the commutatant of $B$, 
is an associative subalgebra of $A$. If $B$ is
also the commutatnt of $B'$, then  $B$ and $B'$ are called
{\em mutual commutants}.

The notion of mutual commutants is an important concept in the
theory of associative algebras, especially in operator algebras.
It has been generalized to the context of groups by Roger Howe \cite{Howe},
called dual pairs of groups, in his study of representation
theory and mathematical physics.
On the classical level, or  more
precisely on the level of Poisson manifolds, an analogue was introduced by
Weinstein \cite{Weinstein:1983}, which is called dual 
pair of Poisson manifolds.
In fact,  Poisson manifolds $\frakg^*$ and $M/G$, 
together with the  Poisson  maps $J: M\lon \frakg^*$ and $p: M\lon M/G$,
 consist a dual pair in terms of Weinstein
\cite{Weinstein:1983}. 
Now  $(C^{\infty}(\frakg^{*})[[\hbar ]] , *_{\hbar}) $
 provides a deformation quantization for the Lie-Poisson
structure $\frakg^*$ while $(C^{\infty}(M)^{G}[[\hbar ]], *_{\hbar})$
quantizes the reduced Poisson space  $M/G$.  So  Theorem
\ref{thm:dual-pair} equivalently  says that under
  deformation quantization,
 the classical dual pair $\frakg^*$ and $M/G$
becomes   mutual commutants. For this reason, we shall
call the pair of algebras 
 $(C^{\infty}(\frakg^{*})[[\hbar ]] , *_{\hbar}) $
and  $(C^{\infty}(M)^{G}[[\hbar ]], *_{\hbar})$
a {\em quantum dual pair}.
In general, given two Poisson manifolds
$P_{1}$ and $P_{2}$ and their deformation
quantization $(C^{\infty}(P_{1})[[\hbar ]], *_{\hbar})$ and
$(C^{\infty}(P_{2})[[\hbar ]], *_{\hbar})$, we say that
they consist a {\em  quantum  dual pair} if there is a
symplectic manifold $M$ and a star-product $(C^{\infty}(M)[[\hbar ]],
*_{\hbar})$ on $M$, and algebra morphisms $\rho_{1}: C^{\infty}(P_{1})[[\hbar ]]\lon C^{\infty}(M)[[\hbar ]]$ and
$\rho_{2}: C^{\infty}(P_{2})[[\hbar ]]\lon C^{\infty}(M)[[\hbar ]]$
such that  $\rho_{1} (C^{\infty}(P_{1})[[\hbar ]])$
and $ \rho_{2} (C^{\infty}(P_{2})[[\hbar ]])$ are mutual commutants.


Unfortunately, at the moment we   only   know a   few 
examples of    deformation quantizable Poisson manifolds.
In fact,  we do not know any other  examples of quantum dual pairs
besides  this and the   trivial ones.
In particular, it is not clear in general whether
a classical dual pair (or Morita equivalent Poisson
manifolds \cite{Xu:1991}) can  be quantized to a quantum
dual pair or not.  The answer  to all these
questions relies upon  how successful it is  the
 deformation quantization theory   of Poisson manifolds.

\section{Appendix A (Derivations of $*$-algebras)}

In this Appendix, we  collect  some  basic
facts concerning  derivations of a $*$-algebra $C^{\infty}(M)[[\hbar ]]$  on
a symplectic manifold $M$.
Some of them  are well known. However, as we
shall see,   Fedosov method even sheds   new light  on  understanding
these  results.

\begin{defi}
A derivation  of
a $*$-algebra $(C^{\infty}(M)[[\hbar  ]] , *_{\hbar}) $ is a
 formal power series of $\hbar$ with coefficients
being  linear operators on $C^{\infty}(M)$: $\delta=D_{0}+\hbar D_{1}
+\cdots +\hbar^{i}D_{i}+\cdots   $ such that
\begin{equation}
\label{eq:der}
\delta (f*_{\hbar}g)=\delta f*_{\hbar}g+f *_{\hbar} \delta g , \forall   f, g\in C^{\infty}(M)[[\hbar  ]].
\end{equation}

A derivation is said to be inner if $\delta =ad \ih H=[\ih H , \cdot ]$ for some $H\in C^{\infty}(M)[[\hbar ]]$.
\end{defi}

Suppose that $\delta=\sum_{i} \hbar^{i}D_{i}$
is a derivation.
By expanding both sides of Equation (\ref{eq:der}),
the $\hbar^{0}$ terms yield that
$D_{0}(fg)=D_{0}(f)g+fD_{0}(g)$. That  is, $D_{0}$ is
a vector field.
By considering the $\hbar^{1}$ terms, one obtains that
$$D_{0}\{f, g\}-\{D_{0}f, g\}-\{f, D_{0}g\}=
(D_{1}f)g+f(D_{1}g)-D_{1}(fg). $$
Since the LHS is skew-symmetric with respect to $f$ and $g$
while the RHS is symmetric, both terms have to vanish identically.
Thus, $D_{1}$ is a vector field and $D_{0}$ is a
symplectic vector field.

Moreover, we have the following (see \cite{BFFLS} for an equivalent 
result, which  is  however  concerning    derivations of the corresponding
 deformed Lie algebra; see also \cite{BCG}).

\begin{thm}
\label{thm:derivation}
Suppose that $D=D_{0}+ \hbar
D_{1} +\cdots \hbar^{i}D_{i}+\cdots  $ is a derivation of  a $*$-algebra
 $(C^{\infty}(M)[[\hbar ]] , *_{\hbar})$  on a
symplectic manifold $M$. Then,
\begin{enumerate}
\item the operator  $D_{i}$, for each $i$,
 is a differential operator, and in particular
$D_{0}$ is a symplectic vector field;
\item there is a canonical one-one correspondence
between derivations and $Z^{1}(M)[[\hbar ]]$;
\item under such a  correspondence, inner
derivations correspond to  exact  1-forms in $B^{1}(M)[[\hbar ]]$.
\end{enumerate}
\end{thm}

To begin with, we will consider Fedosov algebras, and assume that $A=W_{D}$ for
some Fedosov connection $D$. We need a couple of lemmas first.

\begin{lem}
\label{lem:dk}
Let $K\in \gm W$ be a section. Then,
$\rho =ad\ih K=[\ih K, \cdot ]$ defines a derivation of $W_{D}$
iff $DK$ is  a scalar closed one-form on $M$.
\end{lem}
\pf It is clear that $\rho $, defined in this way,
 satisfies
the derivation property.

For any $a\in W_{D}$,
$$D\rho (a)=[\ih DK , a]+[\ih K, Da ].$$
If $\rho W_{D}\subset W_{D}$, it follows  that $[DK, a]=0, \ \forall a\in
W_{D}$. Thus $DK=\theta $ is a scalar one-form on $M$,
i.e., $\theta \in \Omega^{1}(M)[[\hbar ]]$. Then $\theta $ must
be closed since $d\theta =D\theta =D^{2}K=0$.

The converse follows essentially from the same argument backwards. \qed

\begin{lem}
\label{lem:der-extension}
For   a $*$-product on a symplectic manifold $M$,
any symplectic vector field $X$  may  extend 
 to a derivation $\delta =\sum_{i} \hbar^{2 i}D_{2   i}$,
with $D_{2i}$ being differential operators,  such that $D_{0}=X$.
 If $X$ is a hamiltonian
vector field, $\delta$ may  be chosen as an inner derivation.
\end{lem}
\pf Since every $*$-product is equivalent to a Fedosov $*$-product,
we may confine ourselves to a Fedosov algebra $W_{D}$.

Let $\theta =X\per \omega $. Then $\theta $ is a
closed one-form on $M$. Let $K\in \gm (W)$ be
the section satisfying
\begin{equation}
\label{eq:k0}
DK=\theta  \ \ \mbox{and } K|_{y=0}=0.
\end{equation}
 Note that such a section always
exists according to the  Fedosov iteration method.
In fact, $K$ is uniquely determined by
the following  iteration formula:
\begin{equation}
K=-\delta^{-1}\theta +\delta^{-1}(\partial K+[\ih \gamma , K]) .
\end{equation}
Take $\delta  =ad\ih K $. Then $\delta$  is  easily seen to be  a  required
derivation. In fact, locally $\theta $ is 
exact, so $\delta$ can be expressed as an inner derivation generated by  some function
on $M$, which is clearly a  formal power series
of $\hbar^2$ with coefficients being differential operators.

If $X$ is a hamiltonian vector field  with the hamiltonian function $H$,
we may take  $\delta =[\ih H, \cdot ]$, which  is an inner derivation having the
desired property. \qed
{\bf Proof of  Theorem \ref{thm:derivation}}\ \ According to
the observation preceding Theorem \ref{thm:derivation}, $D_{0}$ is a
symplectic vector field. Thus,
 it extends to a derivation $\delta_{0}=D_{0}+O(\hbar )$, whose
coefficients are differential operators,  according to the lemma
above.
Let  $  \tilde{\delta }=\frac{1}{\hbar} (\delta -\delta_{0})$.
Then $  \tilde{\delta } $ is  a derivation and
 $\delta =\delta_{0} +\hbar \tilde{\delta } $.
Applying this process  repeatedly, one obtains
that $\delta =\sum \hbar^{i}\delta_{i}$, where every
$\delta_{i}$ is  a  derivation  whose coefficients
are differential operators. So $\delta$ itself is  such a
derivation as well.

 To continue, without loss of generality,
 we  shall confine ourselves to the case of a Fedosov
 algebra $W_{D}$. Assume that $\delta : W_{D}\lon W_{D}$
is a derivation. For any $x_{0}\in M$, we define a derivation
 $\rho_{x_{0}}$ on $W_{x_{0}}$ by

\begin{equation}
\rho_{x_{0}} a(y ,\hbar )=(\delta \tilde{a})|_{x=x_{0}}, \ \ \ \ \mbox{where }
a(y ,\hbar ) \in W_{x_{0}} \ \mbox{and }
\tilde{a} \in W_{D} \ \mbox{ such that }
\tilde{a}(x_{0}, y, \hbar )=a(y, \hbar ).
\end{equation}
Clearly,
$\rho_{x_{0}} $ is well defined since  $\delta $ is a local operator according
to Part (i).
Since all derivations on $W_{x_{0}}$ are inner, there is an
element  $K(x_{0})\in W_{x_{0}}$ such that
$\rho_{x_{0}}=ad\ih K(x_{0})$. By requiring
that $K(x_{0})|_{y=0}=0$,  $K(x_{0})$ will be unique.
Repeating this process pointwisely, we obtain a
 global section $K\in \gm (W)$ with $K|_{y=0}= 0$ such that

$$\delta \tilde{a}=[\ih K , \tilde{a} ], \ \ \  \ \forall \tilde{a}\in W_{D}. $$
According to Lemma \ref{lem:dk}, $\theta =DK$ belongs to
$Z^{1}(M)[[\hbar ]]$. In this way, we obtain a map $\phi $
 from the  space of derivations  $Der C^{\infty}(M)[[\hbar ]]$ to
that of closed 1-forms  $Z^{1}(M)[[\hbar ]]$.

Conversely, given any $\theta \in  Z^{1}(M)[[\hbar ]]$,    the equation:

\begin{equation}
\label{eq:k}
DK=\theta , \ \ \ \mbox{ and } K|_{y=0}=0
\end{equation}
always has  a unique solution. Thus,
$\delta =ad \ih K$ defines a derivation of $W_{D}$ such that $\phi \delta =\theta $.
In other words, $\phi$ is onto.
It is also simple to see that $\phi$ is injective since the
solution to Equation (\ref{eq:k}) is unique.

Suppose that $\delta$ is an inner derivation:
 $\delta =ad \ih H$ for some $H\in W_{D}$. Thus, $K=H-H_{0}$,  where $H_{0}=H|_{y=0}\in
C^{\infty}(M)[[\hbar ]]$, and
$\theta =DK=D(H-H_{0})=-dH_{0}$, which is clearly exact.

Conversely, if $ \phi \delta =\theta $ is  exact, i.e., $\theta =dH$ for
some $H\in C^{\infty}(M)[[\hbar ]]$,  we have
$D(K-H)= \theta -dH=0$. That is, $K-H \in W_{D}$.
Thus, $\delta =ad  \ih K=ad \ih (K-H)$ is clearly inner.
This concludes the proof of the theorem. \qed

It is well known that the bracket of any two  symplectic vector
fields is hamiltonian.
As  another immediate  application of 
Fedosov method, we obtain the  following ``quantum" analogue of this fact.

\begin{pro}
\label{pro:exact}
Let $*_{\hbar}$ be any $*$-product on a symplectic
manifold $M$. Then the bracket of any derivations
is an inner derivation.
\end{pro}
\pf  Without loss of generality, assume that $A=W_{D}$ for
some Fedosov connections, and
$\delta_{1}=ad{\ih  K_{1}} $ and $\delta_{2}=ad{\ih  K_{2}}$ are
derivations of $W_{D}$,
where $K_{1}$ and $K_{2}$ are sections of $W$. Then, $DK_{1}$ and $DK_{2}$
belong to $Z^{1}(M)[[\hbar ]]$ according to  Lemma \ref{lem:dk}.
Now $[\delta_{1} , \delta_{2}]=ad  \ih (\ih [K_{1} ,  K_{2}])$.
It is clear that $K=\ih [K_{1} ,  K_{2}]$ is a section of $W$ and
$DK=0$. Therefore, $[\delta_{1} , \delta_{2}]$ is an inner
derivation. \qed
{\bf Remark} (1). We note that our definition of inner derivations differs
from the usual one by a factor $\ih$.  If we modify the
usual Hochschild coboundary operator by multiplying
the factor $\ih$, Theorem  \ref{thm:derivation} implies that the first
Hochschild cohomology $H^{1}(C^{\infty}(M)[[\hbar ]], C^{\infty}(M)[[\hbar ]])$ is
 isomorphic to $H^{1}(M)[[\hbar ]]$ (see \cite{BFFLS}).
Similarly, for  higher order cohomology, we expect that
$H^{*} (C^{\infty}(M)[[\hbar ]], C^{\infty}(M)[[\hbar ]]) \cong H^{*} (M) [[\hbar ]]$.
Under such an isomorphism, the Gerstenhaber bracket on  Hochschild cohomology
 would go   to zero on the right hand side.  Moreover, 
the isomorphism at  the second order should
 provide an intrinsic  explanation
for   the characteristic class of a  $*$-algebra  (\cite{WeinsteinX:1996}).

(2). Using  the identification as in  Theorem \ref{thm:derivation}, one obtains
a Lie  bracket on the space $Z^{1}(M)[[\hbar ]]$ such that the bracket
of any closed one-forms is exact (see  Proposition \ref{pro:exact}).
It is easy to see that for any $\theta_{1} ,\ \theta_{2}$,
 $[\theta_{1} , \theta_{2}]=\{\theta_{1} , \theta_{2}\}+O(\hbar )$,
where $\{ \cdot , \cdot \}$ refers to the standard Lie bracket
on one-forms induced from the Poisson bracket (see \cite{AM}). However, it is
difficult to find an explicit expression  for the entire bracket
$[\cdot , \cdot]$. The latter   should be related to
the Weyl curvature of the deformation. 

 Also, for a symplectic manifold, it is well known that
the Poisson  bracket defined on closed one-forms extends to a bracket
on  all one-forms.  It is not clear, however,  whether one can extend
  the bracket $[\cdot , \cdot]$ above  to $\Omega^{1}(M)[[\hbar ]]$. 
It seems that these problems  are all 
 related to the  question raised  by Weinstein  regarding
``quantum Lie algebroids" \cite{Weinstein:1996}.

Another interesting consequence  of Theorem \ref{thm:derivation} is the following:

  \begin{cor}
If $\delta =D_{0}+\hbar D_{1}+\hbar^{2} D_{2} +\cdots =\sum_{i}
\hbar^{i}D_{i}   $ is
a derivation. Then both
$\delta_{even}=D_{0}+\hbar^{2} D_{2} +\cdots =\sum_{i}\hbar^{2i}D_{2i}$
and $\delta_{odd}=D_{1}+\hbar^{2}D_{3} +\cdots =\sum_{i}\hbar^{2i}D_{2i+1}$,
are derivations.
\end{cor}
\pf Without loss of generality,  we  assume that this is a Fedosov $*$-product, and
consider $W_{D}$ as our algebra.

Then each  derivation $\delta $ can be written as
 $\delta =[\ih K, \cdot ] $,  for some $ K\in  \gm W$. Let $\theta =DK \in
Z^{1}(M)[[\hbar ]]$.   Write $\theta =\theta_{even}+\hbar \theta_{odd}$, where
$\theta_{even}$ is the sum of all even terms in $\hbar$ while
$\theta_{odd}$ is the sum of all odd terms in $\hbar$ divided by $\hbar$.
Let $K_{even }$ and $K_{odd}$ be the sections of $W$ corresponding to
$\theta_{even}$ and $\theta_{odd}$,  respectively,  given by Equation (\ref{eq:k}), and
let
\be
\tilde{ \delta}_{even}&=&ad \ih K_{even } \ \ \mbox{and}\\
\tilde{ \delta}_{odd}&=&ad \ih K_{odd} .
\ee

Then,  clearly both $\tilde{ \delta}_{even}$ and $\tilde{ \delta}_{odd}$
 are derivations, and consist only even powers of $\hbar$ according to
Lemma \ref{lem:der-extension}.
Also  $\delta = \tilde{ \delta}_{even}+\hbar  \tilde{ \delta}_{odd}$ by construction.
Therefore, $\tilde{ \delta}_{even}=\delta_{even}$ and
$\tilde{ \delta}_{odd}=\delta_{odd}$.
This concludes the proof. \qed

\section{Appendix B}
A $*$-product
$$f*_{\hbar}g=\sum_{k}\hbar^{k}C_{k}(f, g)$$
on a regular Poisson manifold is said to be 
{\em \preferred } if $C_{k} (\cdot , \cdot )$  is a
 leafwise  bidifferential
operator for every $k$.

In this section, for  completeness, 
we will outline a proof for the following result,
whose proof 
  can also  be found in \cite{NT1}.

\begin{pro}
\label{prop:main}
Suppose that $P$ is a regular Poisson manifold whose symplectic
foliation is  a fibration $P\lon M$.
Assume that the second  Betti number of the   fibers is zero. Then there exists
essentially a unique \preferred  $*$-product on $P$.
I.e., any two \preferred  $*$-products
are equivalent,  and the equivalence  operator
can be chosen as a  formal power series in $\hbar$ with
coefficients being   leafwise differential operators.
\end{pro}

Recall that,   given a fibration $\calf$ on a manifold $P$,
 a leafwise Hochschild cochain is  a $k$-linear
form on $C^{\infty}(P)$ with value in $C^{\infty}(P)$, which
requires to be a  leafwise $k$-differential operator on $P$.
The Hochschild differential is given by
\be
(\delta c)(u_{0}, \cdots , u_{k})&=&u_{0}c(u_{1}, \cdots , u_{k})+
\sum_{i=0}^{k-1}(-1)^{i+1}c(u_{0}, \cdots , u_{i}u_{i+1}, \cdots , u_{k})\\
&&+(-1)^{k+1}c(u_{0}, \cdots , u_{k-1})u_{k} .
\ee

As usual, the space of $k$-cochains is denoted by
$C^{k}_{\calf}(C^{\infty}(P), C^{\infty}(P))$, while
its cohomology is denoted by
$H^{k}_{\calf}(C^{\infty}(P), C^{\infty}(P))$.

\begin{lem}
(Nest-Tsygan \cite{NT1})

Let $\calf$ be a fibration on  a manifold $P$. Then,
\begin{equation}
 H^{k}_{\calf}(C^{\infty}(P), C^{\infty}(P)) \cong \gm  (\wedge^{k}T\calf).
\end{equation}

\end{lem}

An immediate consequence is the following

\begin{lem}
\label{lem:V}
\begin{enumerate}
\item If $c \in C^{2}_{\calf}(C^{\infty}(P), C^{\infty}(P))$ is an antisymmetric two-cocycle, then $c$
is a leafwise bivector field, i.e., $c\in \gm (\wedge^{2}T\calf )$.
\item If $c\in C^{2}_{\calf}(C^{\infty}(P), C^{\infty}(P))$ is
a  symmetric  two-cocycle, then it is a coboundary.
\end{enumerate}
\end{lem}
{\bf Proof of Proposition \ref{prop:main}}
The proof is simply  a  modification of the standard  one for symplectic
manifolds. The idea, roughly speaking,  is as follows.
The classification of $*$-products  is equivalent to that
of their  commutator Lie algebras $[\cdot , \cdot ]_{*}$,
which are  Lie algebra deformations  of the Poisson bracket
 $\{\cdot , \cdot \}$.
The latter is classified by the 2nd order leafwise
Chevalley cohomology of the Poisson Lie algebra.
However, the 2nd order leafwise Chevalley cohomology is isomorphic
to the second leafwise de-Rham cohomology,  and is zero by  assumption.

Recall that  the  Chevalley coboundary (see \cite{L})  operator is
given  by:

\begin{equation}
(\partial c)(u_{0}, \cdots , u_{p} )=\epsilon_{0\cdots  p}^{\lambda_{0}
\cdots  \lambda_{p}}[\frac{1}{p!}\{u_{\lambda_{0}}, c(u_{\lambda_{1}},
\cdots , u_{\lambda_{p}}) \}-\frac{1}{2(p-1)!}c(\{u_{\lambda_{0}}, u_{\lambda_{1}}\},
u_{\lambda_{2}}, \cdots , u_{\lambda_{p}}) ],
\end{equation}
where $u_{\lambda} \in C^{\infty}(P)$ and $\epsilon$ is  the Kronecker symbol.
If $c\in C^{2}_{\calf}(C^{\infty}(P), C^{\infty}(P))$ is antisymmetric,
$$(\partial c)(u_{0}, u_{1}, u_{2})=
(\{u_{0}, c(u_{1}, u_{2})\}-c(\{u_{0}, u_{1}\}, u_{2}))+c.p.,
$$
where $c.p.$ stands for the cyclic permutation.
Hence if $c\in \gm(\wedge^{2} T \calf )$,
$\partial c=[\pi , c]=d_{\pi }c$,
where the bracket is the Schouten bracket on multivector fields and
$d_{\pi }: \gm (\wedge^{*}T \calf ) \lon \gm (\wedge^{*+1}T \calf )$ is
the differential operator defining the (leafwise)-Poisson cohomology.

Suppose that
\be
f*g&=& \sum \hbar^{k} C_{k}(f, g), \ \ \mbox{and }\\
f*'g&=& \sum \hbar^{k} C_{k}'(f, g)
\ee
are two \preferred $*$-products on $P$.

We need to construct an equivalence operator between them.
This will   proceed   by induction.
Assume that they are equivalent up to $i=2k$, i.e., we can find
an equivalence operator  under which, $C_{i}=C_{i}'$ for $0\leq i\leq 2k$.
Then, $\delta (C_{2k+1}-C'_{2k+1})=0$. Since
$C_{2k+1}-C'_{2k+1}$ is skew-symmetric,
it belongs to $ \gm (\wedge^{2}T\calf )$ by Lemma \ref{lem:V}.

On the other hand, both $[\cdot , \cdot ]_{*}$
and $[\cdot , \cdot ]_{*'}$ are deformations
of the Poisson  Lie algbra $\{\cdot , \cdot \}$.
 Thus, $\partial (C_{2k+1}-C'_{2k+1})=0$.
That is, $[\pi , C_{2k+1}-C'_{2k+1} ]=0$. In other words, $C_{2k+1}-C'_{2k+1}$
is a 2-cocycle of the (leafwise)-Poisson cohomology, which is
isomorphic to the leafwise de-Rham cohomology. Since the 2nd leafwise de-Rham cohomology
is zero by assumption,
  there is  a vector field  $X\in \gm (T\calf )$ such that
$  C_{2k+1}-C'_{2k+1}=[\pi , X]$. In other words,
$$C_{2k+1}(f, g)-C'_{2k+1}(f, g)=\{Xf ,g\}+\{f, Xg\}-X\{f ,g\}, \ \ \forall
f , g\in C^{\infty}(P). $$
It is then easy to see that $T=1+\hbar^{2k}X$ establishes
an isomorphism between $*$ and $*'$ up to $\hbar^{2k+1}$.
That is,
$$T(f*g)-Tf*'Tg=O(\hbar^{2k+2}). $$
Finally, we assume  that $*$ and $*'$ are equivalent  up to $n=2k-1$.
By applying an equivalence operator, we may assume that
$C_{i}=C_{i}'$ for $0\leq i\leq 2k-1$.
Then, $\delta (C_{2k}-C'_{2k})=0$. Since $C_{2k}-C'_{2k}$
is a symmetric 2-cochain, $c_{2k}-c'_{2k}=\delta D$,
for some  $D\in C^{1}_{\calf}(C^{\infty}(P), C^{\infty}(P))$ according
to Lemma \ref{lem:V}.
Thus $T=1+\hbar^{2k}D$ establishes an isomorphism between
$*$ and $*'$ up to $\hbar^{2k}$. This concludes the proof. \qed

\section{Appendix C}

Recall that a {\em \nat } $*$-product  on a symplectic manifold $M$ is
a $*$-product:
$$ u*_{\hbar }v=uv-\frac{i\hbar }{2}\{u , v\} + \half (-\frac{i\hbar }{2})^{2}Q_{2}(u, v)+\cdots , $$
where $Q_{2}(u, v)$ is a bidifferential operator of order $2$ in each argument.

It is well known \cite{BFFLS} \cite{L} that associated to a  \nat 
$*$-product there
is a canonical torsion-free symplectic connection. More precisely, we
have

\begin{pro}
\label{pro:nat-split}
Let $ u*_{\hbar }v=uv-\frac{i\hbar }{2}\{u , v\} +\half (-\frac{i\hbar }{2})^{2}Q_{2}(u, v)+\cdots $
be a \nat $*$-product on $M$. Then there exists a unique torsion-free symplectic connection
$\nabla$ such that
$$Q_{2}(u, v)=P_{\nabla}^{2}(u, v)+H(u, v), $$
where $H(u, v)$ is a bidifferential operator of maximum order  $1$ in each argument.
\end{pro}

For completeness, we  outline a proof here  in this appendix. \\\\\\
{ \bf Proof of Proposition \ref{pro:nat-split}} Take a torsion free symplectic
connection $\tilde{\nabla}$. Then,
$$ u*_{\hbar }v=uv-\frac{i\hbar }{2}\{u , v\} +\half (-\frac{i\hbar }{2})^{2} P_{\tilde{\nabla}}^{2}(u, v)$$
 is a deformation of order up to $\hbar^2$. Therefore,
$\delta (Q_{2}- P_{\tilde{\nabla}}^{2})=0$, where $\delta$ denotes the usual
Hochschild differential. Since $Q_{2}- P_{\tilde{\nabla}}^{2}$ is symmetric, it
must be a 2-coboundary. Hence there is a differential operator $D$ of  maximum
order $3$ such that
$$Q_{2}- P_{\tilde{\nabla}}^{2}=\delta D. $$
The principal term of $D$ corresponds to a $3$-covariant symmetric tensor  $T^{ijk}$.

In local coordinates, write $\tilde{\nabla}_{\parr{i}}\parr{j}=\tilde{\gm}_{ij}^{k}
\parr{k}$, and let $\tilde{\gm}^{ijk}=\tilde{\gm}_{lm}^{i}\pi^{lj}\pi^{mk}$.
Set
\begin{equation}
{\gm}^{ijk}=\tilde{\gm}^{ijk}+3 T^{ijk} .
\end{equation}
Since $T^{ijk}$ is a completely symmetric tensor, the equation above
defines a torsion free symplectic connection $\nabla_{\parr{i}}\parr{j}={\gm}_{ij}^{k}\parr{k}$
with $\gm^{ijk}={\gm}_{lm}^{i}\pi^{lj}\pi^{mk}$.
A simple calculation  yields that
\begin{eqnarray}
P_{\tilde{\nabla}}^{2}(u, v)-P_{\nabla}^{2}(u, v)&=&\pi^{i_{1}j_{1}}\pi^{i_{2}j_{2}}
(\tilde{\gm}_{i_{1}i_{2}}^{k}\tilde{\gm}_{j_{1}j_{2}}^{l} -
{\gm}_{i_{1}i_{2}}^{k}{\gm}_{j_{1}j_{2}}^{l} )\parrr{u}{k}\parrr{v}{l}  \nonumber \\
&& + 3T^{i_{1}i_{2}i_{3}}\parrr{u}{i_{1}}\frac{\partial^{2}v}{\partial x^{i_{2}} \partial x^{i_{3}}}
+3T^{i_{1}i_{2}i_{3}} \frac{\partial^{2}u}{\partial x^{i_{1}} \partial x^{i_{2}}} \parrr{v}{i_{3}}.
\label{eq:difference}
\end{eqnarray}

On the other hand,
$$(\delta D)(u, v)=-3T^{i_{1}i_{2}i_{3}}\parrr{u}{i_{1}}\frac{\partial^{2}v}{\partial x^{i_{2}}
 \partial x^{i_{3}}}-3T^{i_{1}i_{2}i_{3}} \frac{\partial^{2}u}{\partial x^{i_{1}} \partial x^{i_{2}}} \parrr{v}{i_{3}}+
\tilde{H}(u, v),$$
where $\tilde{H}(u, v)$ is    a bidifferential operator of maximum order  $1$ in each argument.
Therefore,
$$  Q_{2}(u, v)-P_{\nabla}^{2}(u, v)=P_{\tilde{\nabla}}^{2}(u, v)-P_{\nabla}^{2}(u, v)+(\delta D)(u, v)$$
is clearly a  bidifferential operator of maximum order  $1$ in each argument.

To see that such a connection is unique, it
suffices to note that $P_{\tilde{\nabla}}^{2}(u, v)-P_{\nabla}^{2}(u, v)$
 is a   bidifferential operator of maximum order  $1$  iff
$\tilde{\nabla}=\nabla$. This can be  easily seen from Equation (\ref{eq:difference}). \qed

\end{document}